\documentclass[10pt,a4paper,reqno]{amsart}
\usepackage{tabls}
\usepackage{url}
\usepackage[linktocpage = true]{hyperref}
\hypersetup{
 colorlinks = true,
 citecolor = blue,
}
\usepackage{color}
\definecolor{red}{rgb}{1,0,0}
\definecolor{green}{rgb}{0,1,0}
\definecolor{blue}{rgb}{0,0,1}
\definecolor{refkey}{gray}{.625}
\definecolor{labelkey}{gray}{.625}
\usepackage[lite]{amsrefs}
\usepackage[a4paper]{geometry}
\usepackage{amsfonts}
\usepackage{amssymb}
\usepackage{mathrsfs}
\usepackage{amsmath}
\usepackage{amsthm}
\usepackage{amscd}
\usepackage{enumerate}
\usepackage{paralist}
\usepackage{xypic}
\usepackage{graphicx}
\usepackage{mathtools}
\usepackage[all,cmtip]{xy}
\usepackage{kantlipsum}
\allowdisplaybreaks

\newcommand{\abs}[1]{\lvert#1\rvert}

\newcommand{\At}{\operatorname{At}}
\newcommand{\A}{\mathscr{A}}
\newcommand{\B}{\mathcal{B}}

\newcommand{\id}{\operatorname{id}}

\newcommand{\Z}{\mathbb{Z}}

\newcommand{\g}{\mathfrak{g}}

\newcommand{\bd}{\begin{displaymath}}
\newcommand{\ed}{\end{displaymath}}
\newcommand{\be}{\begin{equation}}
\newcommand{\ee}{\end{equation}}

\newcommand{\DG}{\mathrm{dg}}

\makeatletter
 \def\title@font{\normalsize\bfseries}
 \let\ltx@maketitle\@maketitle
 \def\@maketitle{\bgroup%
 \let\ltx@title\@title%
 \def\@\title{\resizebox{\textwidth}{!}{%
  \mbox{\title@font\ltx@title}%
 }}%
 \ltx@maketitle%
 \egroup}
\makeatother

\theoremstyle{plain}
\newtheorem*{zorn*}{Zorn's lemma}
\newtheorem*{tychonoff*}{Tychonoff's theorem}

\newtheorem{lem}[equation]{Lemma}
\newtheorem{Cor}[equation]{Corollary}
\newtheorem{Thm}[equation]{Theorem}
\newtheorem*{theorem*}{Theorem}
\newtheorem{Def}[equation]{Definition}

\newtheorem{def-prop}[equation]{Definition-Proposition}
\newtheorem{prop}[equation]{Proposition}
\newtheorem{prop-def}[equation]{Proposition-Definition}
\newtheorem{Ex}[equation]{Example}
\newtheorem{Rem}[equation]{Remark}
\newtheorem{Convention}[equation]{Convention}

\numberwithin{equation}{section}

\begin{document}
\def\bp{\begin{proof}}
\def\ep{\end{proof}}
\def\be{\begin{equation}}
\def\ee{\end{equation}}
\def\C{\mathbb{C}}
\def\CE{\mathrm{CE}}
\def\D{\mathcal{D}}
\def\E{\mathscr{E}}
\def\F{\mathscr{F}}
\def\H{\textbf{H}}
\def\k{\mathbb{K}}
\def\G{\mathcal{G}}
\def\M{\mathcal{M}}
\def\m{\mathfrak{m}}
\def\t{\mathfrak{t}}
\def\O{\mathcal{O}}
\def\P{\mathcal{P}}
\def\U{\mathcal{U}}
\def\V{\mathscr{V}}
\def\L{\mathcal{L}}
\def\X{\mathbb{X}}
\def\Y{\mathbb{Y}}
\def\spec{\text{spec}}
\def\Im{\text{Im}}
\def\coker{\operatorname{coker}}
\def\Ext{\operatorname{Ext}}
\def\End{\operatorname{End}}
\def\Kap{\operatorname{Kap}}
\def\pr{\operatorname{pr}}
\def\id{\operatorname{id}}
\def\Der{\operatorname{Der}}
\def\Hom{\operatorname{Hom}}
\def\Jet{\operatorname{Jet}}
\def\Map{\operatorname{Map}}
\def\Mod{\operatorname{Mod}}
\def\sgn{\operatorname{sgn}}
\def\sh{\operatorname{sh}}

\newcommand{\CATderivationA}{\operatorname{dgDer}_\A}

\newcommand{\CATLeibnizoneA}{\operatorname{Leib}_\infty(\A)}

\newcommand{\CAThtyAmod}{\operatorname{H}(\operatorname{dg}\A)}

\title{Kapranov's Construction of sh Leibniz Algebras}

\author{Zhuo Chen}
\address{Department of Mathematics, Tsinghua University} 
\email{\href{mailto:zchen@math.tsinghua.edu.cn}{zchen@math.tsinghua.edu.cn}}
\thanks{Research partially supported by NSFC grant 11471179}

\author{Zhangju Liu}
\address{Department of Mathematics, Peking University} 
\email{\href{mailto:liuzj@pku.edu.cn}{liuzj@pku.edu.cn}}

\author{Maosong Xiang}
\address{Center for Mathematical Sciences, Huazhong University of Science and Technology}
\email{\href{mailto: msxiang@hust.edu.cn}{msxiang@hust.edu.cn}}

\begin{abstract}
Motivated by Kapranov's discovery of an sh Lie algebra structure on the tangent complex of a K\"{a}hler manifold and Chen-Sti\'{e}non-Xu's construction of sh Leibniz algebras associated with a Lie pair, we find a general method to construct sh Leibniz algebras. Let $\A$ be a commutative dg algebra. Given a derivation of $\A$ valued in a dg module $\Omega$, we show that there exist sh Leibniz algebra structures on the dual module of $\Omega$. Moreover, we prove that this process establishes a functor from the category of dg module valued derivations to the category of sh Leibniz algebras over $\A$.   
\emph{MSC class}:~16E45, 18G55.
\end{abstract}

\maketitle

\tableofcontents
\section{Introduction}\label{Sec:First}

Higher homotopies and higher structures are playing  important roles in mathematics and  some branches of theoretical physics, such as gauge theory and topological field theory (see Huebschmann~\cite{Huesurvey}). Higher homotopies, as explained by Huebschmann in~\cite{HueshLie}, often arise from the process of transferring certain strict geometric or algebraic structure on a huge chain complex to a smaller but chain homotopic complex. For instance, an sh Lie algebra (also known as $L_\infty$-algebra~\cite{LS}) yields  from a dg Lie algebra by applying homological perturbation theory~\cite{HueLie}. Here and in the sequel, sh is short for \emph{strongly homotopy} and dg is short for \emph{differential graded}.
Sh Leibniz algebras, also known as sh Loday algebras, Loday infinity algebras or Leibniz$_\infty$ algebras~\cite{AP}, are also examples of higher structures. In fact, the notion of   Leibniz$_\infty$ algebras is a generalization of   $L_\infty$ algebras where the skew-symmetricity constraint on multibrackets is discarded.

In this note, we use the notion of Leibniz$_\infty[1]$ algebras (see Definition~\ref{Def:Leibnizinfinity}), which is equivalent to the notion of sh Leibniz algebras, and
study a particular method to construct Leibniz$_\infty[1]$ algebras.
This method first appeared in Kapranov's approach to Rozansky-Witten theory~\cite{Kap}: Given a K\"{a}hler manifold $X$, Kapranov discovered an $L_\infty$ algebra structure on $\Omega^{0,\bullet-1}_X(T_X)$ via the Atiyah class $\alpha_X$. More precisely,
let $\nabla$ be the Chern connection on the holomorphic tangent bundle $T_X$. Then the curvature $R_\nabla \in \Omega_X^{0,1}(\Hom(S^2(T_X), T_X))$ is a Dolbeault representative of the Atiyah class $\alpha_X$. The $L_\infty$ brackets $\{\lambda_k\}_{k\geq1}$ on $\Omega^{0,\bullet-1}_X(T_X)$ are defined by
\begin{itemize}
\item  $\lambda_1 = \bar{\partial}$.
\item  $\lambda_2 = R_\nabla$.
\item  $\lambda_{k+1} = \nabla^{1,0}(\lambda_k) \in \Omega_X^{0,1}(\Hom(S^{k+1}(T_X), T_X))$, {for } $k = 2,3,4,...$.
\end{itemize}

 Kapranov's construction of $L_\infty$ algebras is generalized in Chen, Sti\'{e}non and Xu's work~\cite{CSX} where the setting is a Lie algebroid pair (Lie pair, for short) $(L,A)$. It is shown that the graded  vector space $\Gamma(\wedge^\bullet A^\vee \otimes L/A)$ admits a Leibniz$_\infty[1]$ algebra structure (\cite{CSX}*{Theorem 3.13}) via the Atiyah class of the Lie pair $(L,A)$. The Atiyah class of Lie pairs encompasses the original Atiyah class~\cite{Atiyah} of holomorphic vector bundles and the Molino class~\cite{Molino1} of foliations as special cases.
This construction of Leibniz$_\infty[1]$ algebra structures is similar to that of Kapranov
 --- First, we choose a splitting $j: L/A \to L$ of vector bundles so that $L \cong A\oplus L/A$.
 Second, choose an $L$-connection $\nabla$ on $L/A$ extending the $A$-module structure.
 Then the Leibniz$_\infty[1]$ brackets $\{\lambda_k\}_{k\geq1}$ on $\Gamma(\wedge^\bullet A^\vee \otimes L/A)$ are determined as follows:
\begin{itemize}
  \item $\lambda_1 = d_{\CE}$ is the Chevalley-Eilenberg differential of the Bott representation of $A$ on $L/A$.
  \item Define a bundle map $R_2: L/A \otimes L/A \rightarrow A^\vee \otimes L/A$ via the Atiyah cocycle $\alpha_{L/A}^{\nabla}$ (see Section \ref{Ex: Lie pair}):
      \bd
       R_2(b_1,b_2) = \alpha_{L/A}^{\nabla}(-,b_1)b_2,\;\;\forall b_1,b_2 \in \Gamma(L/A).
      \ed
      The second structure map $\lambda_2$ is specified by
     \bd
      \lambda_2(\xi_1 \otimes b_1, \xi_2 \otimes b_2) = (-1)^{\abs{\xi_1}+\abs{\xi_2}} \xi_1\wedge\xi_2 \wedge R_2(b_1,b_2),
     \ed
     for all $\xi_1,\xi_2 \in \Gamma(\wedge^\bullet A^\vee)$ and $b_1,b_2 \in \Gamma(L/A)$.
  \item Define a sequence of bundle maps $R_k: (L/A)^{\otimes k} \rightarrow A^\vee \otimes L/A, k \geq 3$ recursively by $R_{n+1}=\nabla R_n$, i.e.,
      \bd
       R_{n+1}(b_0 \otimes \cdots \otimes b_n) = R_n(\nabla_{j(b_0)}(b_1 \otimes \cdots \otimes b_n)) - \nabla_{j(b_0)}R_n(b_1 \otimes \cdots \otimes b_n).
      \ed
      The $k$-th structure map is specified by
      \bd
       \lambda_k(\xi_1 \otimes b_1,\cdots,\xi_k\otimes b_k) = (-1)^{\abs{\xi_1}+\cdots+\abs{\xi_k}}\xi_1\wedge\cdots\wedge\xi_k\wedge R_k(b_1,\cdots,b_k),
      \ed
      for all $\xi_i \in \Gamma(\wedge^\bullet A^\vee), b_i \in \Gamma(L/A), 1 \leq i \leq k$.
\end{itemize}
We call $(\Gamma(\wedge^\bullet A^\vee \otimes L/A),\{\lambda_k\}_{k\geq1})$ a Kapranov Leibniz$_\infty[1]$ algebra. Its construction needs, \emph{a priori}, some extra choices (a splitting $j$ and an $L$-connection $\nabla$ on $L/A$). Then one   asks    a natural  question (\cite{CSX}*{Remark 3.19}) --- how does the Leibniz$_\infty[1]$ algebra structure on $\Gamma(\wedge^\bullet A^\vee \otimes L/A)$ depend on the choice of splitting data and connections?
The main goal of this note is to answer this question---Kapranov Leibniz$_\infty[1]$ algebra structures on $\Gamma(\wedge^\bullet A^\vee \otimes L/A)$ associated with different choices of $j$ and $\nabla$, are  mutually isomorphic in the category of Leibniz$_\infty[1]$ algebras over $\Gamma(\wedge^\bullet A^\vee)$ (see Theorem 1.3 or Theorem~\ref{Thm:Liepair}).

We adopt an algebraic approach to achieve this goal.
The algebraic notion we need is a dg module valued derivation of a commutative differential graded algebra (cdga for short) $\A$ (see Definition~\ref{Def:DGderivations}).
As an immediate example from complex geometry, consider a complex manifold $X$. The Dolbeault dg algebra $\A=(\Omega_X^{0,\bullet},\bar{\partial})$ is a cdga. Let $\Omega=(\Omega_X^{0,\bullet}((T^{1,0}X)^\vee), \bar{\partial})$ be the dg $\A$-module generated by the section space of holomorphic cotangent bundle $(T^{1,0}X)^\vee$. Then $\partial: \A \to \Omega$ is an $\Omega$-valued derivation of $\A$.

We now explain how Kapranov's original method and Chen-Sti\'{e}non-Xu's construction can be further generalized in the setting of a dg module valued derivation $\A \xrightarrow{\delta} \Omega$.
Consider the dual dg $\A$-module $\B=\Omega^\vee$ of $\Omega$.
First, one chooses a $\delta$-connection $\nabla$ on $\B$, i.e., a map $\nabla:\B \to \Omega\otimes_\A \B$ that extends the $\delta$-map (see Definition \ref{Def-delta-connection}).
Then one can define a sequence of degree $1$ maps $\{\mathcal{R}^\nabla_k:~\otimes^k\B \to \B\}_{k \geq 1}$ as follows:
\begin{itemize}
\item $\mathcal{R}_1^\nabla = \partial_\A $ is the differential on $\B$.
\item $\mathcal{R}^\nabla_2=\At_\B^\nabla: ~\B \otimes_{\A} \B \rightarrow \B$ is the twisted Atiyah cocycle (see Definition \ref{Def: twisted Atiyah cocycle}).
     \item $\mathcal{R}_k^\nabla$ for $ k \geq 3$ are defined recursively by $\mathcal{R}_k^\nabla =\nabla \mathcal{R}_{k-1}^\nabla$ (see Equation~\eqref{Rnabla}).
\end{itemize}
Our first result is the following
\begin{Thm}\label{Thm1.1}
  When endowed with structure maps $\{\mathcal{R}^\nabla_k\}_{k \geq 1}$, the dg $\A$-module $\B$ becomes a Leibniz$_\infty[1]$ $\A$-algebra.
\end{Thm}
Here by saying that $\B$ is a Leibniz$_\infty[1]$ $\A$-algebra, we mean that its higher structure maps $\{\mathcal{R}^\nabla_k\}_{k\geq 2}$ are all $\A$-multilinear.
We emphasise that the Kapranov Leibniz$_\infty[1]$ algebra   $(\B,\{\mathcal{R}^\nabla_k\}_{k\geq1})$ should be treated as an object in the category of Leibniz$_\infty[1]$  $\A$-algebras.   In fact, if we treat $(\B,\{\mathcal{R}^\nabla_k\}_{k\geq 1})$ merely as a Leibniz$_\infty[1]$ algebra over $\k$, it is always isomorphic to a trivial one (see Remark \ref{Rmk:whyoverA}).
We call $(\B,\{\mathcal{R}^\nabla_k\}_{k\geq1})$ the Kapranov Leibniz$_\infty[1]$ $\A$-algebra associated with the dg module valued derivation $\A \xrightarrow{\delta} \Omega$ and the $\delta$-connection $\nabla$.

Our second result is that the above construction is functorial:
\begin{Thm}\label{Thm: functor}
The above construction defines a functor $\Kap$, called Kapranov functor, from the category of dg module valued derivations of a cdga $\A$ to the category of Leibniz$_\infty$ $\A$-algebras.
Moreover, the Kapranov functor $\Kap$ is homotopy invariant, i.e., if $\delta$ and $\delta^\prime$ are two homotopic derivations of $\A$ valued in the same dg module, then $\Kap(\delta)$ is isomorphic to $\Kap(\delta^\prime)$.
\end{Thm}
Applying Theorem~\ref{Thm: functor} to dg module valued derivations arising from Lie pairs, we obtain the answer of our motivating question:
\begin{Thm}\label{Thm: Lie pair}
  Let $(L,A)$ be a Lie pair. The Kapranov Leibniz$_\infty[1]$ algebra structure on the graded vector space $\Gamma(\wedge^\bullet A^\vee \otimes L/A)$ is unique up to isomorphisms in the category of Leibniz$_\infty[1]$  $\Gamma(\wedge^\bullet A^\vee)$-algebras.
\end{Thm}

This note is organized as follows: Section~\ref{Sec:twistedAtiyahclass} consists of our conventions, notations, and the notion of twisted Atiyah classes. We will see that twisted Atiyah classes encompass  Atiyah classes of Lie pairs and dg Lie algebroids as special cases.
Section~\ref{Sec:MainResults} contains a brief summary of sh Leibniz algebras, the construction of the Kapranov functor, and its applications.
Finally, we present some relevant remarks and open questions in Section~\ref{Sec: final section}.

\textbf{Acknowledgements.}
We would like to thank Bangming Deng, Wei Hong, Kai Jiang, Honglei Lang, Camille Laurent-Gengoux, Mathieu Sti\'{e}non, Jim Stasheff, Yannick Voglaire, and Ping Xu for useful discussions and comments. Xiang is grateful to Penn State University, Peking University and Tsinghua University, for their hospitality during his visits.

\section{Atiyah classes of commutative dg algebras and their twists} \label{Sec:twistedAtiyahclass}

In~\cite{Atiyah}, Atiyah introduced a cohomology class, which has come to be known as Atiyah class, to characterize the obstruction to the existence of holomorphic connections on a holomorphic vector bundle.
The notion of Atiyah classes have been developed in the past decades for diverse purposes (see~\cites{Bottacin,CV,CLX,Costello,CSX,LaurentSX-CR,LaurentSX,MSX}).
In this section, we recall Atiyah classes of commutative dg algebras defined by Costello~\cite{Costello} and introduce a version of twisted Atiyah classes.

\subsection{Atiyah classes of commutative dg algebras}
Throughout this paper, $\k$ denotes a field of characteristic zero and graded means $\Z$-graded.
A commutative differential graded algebra (cdga for short) over $\k$ is a pair $(\A,d_\A)$, where $\A$ is a commutative graded $\k$-algebra, and $d_\A:~\A  \rightarrow \A$, usually called the differential, is a homogeneous degree one derivation of square zero.
We also write $\A$ for a cdga without making its differential explicitly.

An $\A$-module is a representation of the underlying commutative graded algebra of $\A$ by forgetting the differential $d_\A$. 
By a dg $\A$-module, we mean 
an $\A$-module $\E$, together with a degree one and square zero endomorphism $\partial_\A^\E$ of the graded $\k$-vector space $\E$, called the differential, such that
\bd
 \partial^\E_\A  (a e) = (d_\A a )  e + (-1)^{\abs{a}}a \partial_\A^\E(e),
\ed
for all $a \in \A , e \in \E $.
To work with various different dg $\A$-modules, the differential $\partial_\A^\E$ of any dg $\A$-module $\E$ will be denoted by the same notation $\partial_\A$. 
A dg $\A$-module $(\E,\partial_\A)$ will also be simply denoted by $\E$. 

The dg $\A$-module of K\"{a}hler differentials is the graded $\A$-module
$$
\Omega_{\A \mid\k}^1 = \operatorname{span}\{d_{dR}a:~ a \in \A \}/\{d_{dR}(ab)-(d_{dR}a)b - (-1)^{\abs{a}}ad_{dR}b:~ a,b \in \A \},
$$
together with the differential $\partial_\A$ such that the algebraic de Rham operator $d_{dR}:~\A \to \Omega_{\A\mid\k}^1$ is a cochain map, i.e., $\partial_\A (d_{dR}a) = d_{dR}(d_\A a)$ for all $a \in \A$.
In the sequel, we assume that $\Omega^1_{\A\mid\k}$ is projective as an  $\A$-module.

A degree $r$ morphism of dg $\A$-modules, denoted by $\alpha\in \Hom^r_{\DG\A}(\E,\F)$, is a degree $r$ $\A$-module morphism $\alpha:~\E \to \F$, which is also compatible with differentials:
$$
\partial_\A(\alpha) := \partial_\A \circ \alpha -(-1)^{r}\alpha \circ \partial_\A = 0: \E \to \F.
$$

\begin{Def}[Costello~\cite{Costello}]\label{Def Costello}
Let $\A$ be a cdga and $\E$ an $\A$-module.
\begin{compactenum}
  \item A connection on $\E$ is a (degree $0$) map of graded $\k$-vector spaces
  \bd
  \blacktriangledown:~ \E  \rightarrow  \Omega_{\A\mid\k}^1 \otimes_\A \E  ,
  \ed
  satisfying the Leibniz rule
  \bd
   \blacktriangledown(ae) = (d_{dR}a) \otimes e + a\blacktriangledown(e),\;\;\forall a \in \A , e \in \E .
  \ed
  \item Assume that $(\E,\partial_\A)$ is a dg $\A$-module. Given a connection $\blacktriangledown$ on $\E$,
  \bd
   \At_\E^\blacktriangledown :=[\blacktriangledown,\partial_\A] = \blacktriangledown \circ \partial_\A   - \partial_\A \circ \blacktriangledown \in \Omega_{\A\mid\k}^1 \otimes_\A \End_\A(\E)
  \ed
  is a closed element of degree $1$ which measures the failure of $\blacktriangledown$ to be a cochain map. Its cohomology class $\At_\E \in H^1(\A,\Omega_{\A\mid\k}^1 \otimes_\A \End_\A(\E))$ is independent of the choice of connections, and is called the Atiyah class of the dg $\A$-module $\E$.
  \end{compactenum}
\end{Def}
The existence of connections on $\E$ is guaranteed if $\E$ is a projective $\A$-module. Hence we make the following
\begin{Convention}\label{Rmk:projective-exisitence-connections}
 In this note,  all $\A$-modules are assumed to be projective.
\end{Convention}

\begin{Ex}[Mehta-St\'{e}non-Xu~\cite{MSX}]\label{Ex: dg manifold}
  Let $(\M,Q_\M)$ be a smooth dg manifold, where $\M = (M,\O_\M)$ is a smooth $\Z$-graded manifold, and $Q_\M$ is a homological vector field on $\M$. Then $\A = (C^\infty(\M),Q_\M)$ is a cdga. For each dg vector bundle $(\mathcal{E},Q_{\mathcal{E}})$ over $(\M,Q_\M)$, its space of sections $\mathscr{E} = (\Gamma(\mathcal{E}),Q_{\mathcal{E}})$ is a dg $\A$-module. The Atiyah class $\At_\E$ of $\E$ coincides, up to a minus sign, with the Atiyah class $\At_{\mathcal{E}}$ of the dg vector bundle $\mathcal{E}$ with respect to the dg Lie algebroid $T\M$ defined by Mehta-St\'{e}non and Xu. This is a particular instance of Atiyah classes of dg vector bundles with respect to a general dg Lie algebroid (see Section \ref{Sec:DGLieAlgebroidsLiePairs}).
\end{Ex}

\subsection{Dg module valued derivations and twisted Atiyah classes}
A key notion in this note is dg module valued derivation (dg derivation for short):
\begin{Def}\label{Def:DGderivations}
  Let $(\A,d_\A)$ be a cdga and $(\Omega,\partial_\A)$ a dg $\A$-module.
  \begin{itemize}
  \item A dg derivation of $\A$ valued in $(\Omega,\partial_\A)$ is a degree $0$ derivation $\delta:~ \A \rightarrow \Omega$ of the commutative graded algebra $\A$ valued in the $\A$-module $\Omega$,
   \bd
     \delta(ab) = \delta(a)b + a\delta(b),\;\;\;\forall a,b \in \A ,
   \ed
   which commutes with the differentials as well:
  \be\label{compatibility}
     \delta \circ d_\A = \partial_\A \circ \delta:~ \; \A  \rightarrow \Omega.
  \ee
  Such a dg derivation is simply denoted by $\A \xrightarrow{\delta} \Omega$.
  \item Let $\delta$ and $\delta^\prime$ be two $(\Omega,\partial_\A)$-valued dg derivations of $\A$. They are said to be homotopic, written as $\delta\sim \delta'$, if there exists a degree $(-1)$ derivation $h$ of $\A$ valued in the $\A$-module $\Omega$ such that
     \bd
       \delta^\prime - \delta = [\partial_\A,h] = \partial_\A \circ h + h \circ d_\A:  \A \to \Omega.
     \ed
 \end{itemize}
\end{Def}
An immediate example of dg derivations is $ \Omega_X^{0,\bullet}\xrightarrow{\partial} \Omega_X^{0,\bullet} ((T^{1,0}X)^\vee)$ arising from a complex manifold $X$, which has already been explained in Section~\ref{Sec:First}.
Another fundamental example is the dg derivation $\A \xrightarrow{d_{dR}} \Omega_{\A\mid\k}^1$, which is universal in the following sense: For any generic dg derivation $\A \xrightarrow{\delta} \Omega$, there exists a unique dg $\A$-module morphism $\bar{\delta}:~\Omega_{\A\mid\k}^1 \to \Omega$ such that the following diagram commutes:
  \bd
  \xymatrix{
   \A \ar[d]_-{d_{dR}} \ar[r]^-{\delta} & \Omega \\
   \Omega_{\A\mid\k}^1. \ar[ur]_-{\exists ! ~\bar{\delta}} &
  }
  \ed

Thus,
\begin{equation}\label{Eq:bardeltatensor1}
 \bar{\delta} \otimes_\A \id_{\End_\A(\E)}:~ \Omega_{\A\mid\k}^1 \otimes_\A \End_\A(\E) \rightarrow \Omega \otimes_\A \End_\A(\E)
\end{equation}
is a dg $\A$-module morphism as well.
\begin{Def}[Twisted Atiyah class]\label{Def: twisted Atiyah class}
Let $\E$ be a dg $\A$-module and $\A \xrightarrow{\delta} \Omega$ a dg derivation.
The dg $\A$-module morphism $\bar{\delta} \otimes_\A \id_{\End_\A(\E)}$ in Equation~\eqref{Eq:bardeltatensor1} sends the Atiyah class $\At_\E\in H^1(\A,\Omega_{\A\mid\k}^1 \otimes_\A \End_\A(\E))$ of $\E$ to a cohomology class
$$
\At_\E^\delta \in H^1(\A, \Omega \otimes_\A \End_\A(\E)),
$$
which is called the $\delta$-twisted Atiyah class of $\E$.
\end{Def}

It follows immediately that twisted Atiyah classes are homotopic invariant:
\begin{prop}\label{prop: homotopy invariance}
  If $\delta\sim\delta^\prime$, then for any dg $\A$-module $\E$,
  \bd
   \At_\E^\delta = \At_\E^{\delta^\prime} \in H^1(\A, \Omega \otimes_\A \End_\A(\E)).
  \ed
\end{prop}

Below we give a different characterization of the twisted Atiyah class $\At_\E^\delta$. We need another key notion in this note --- $\delta$-connections, which can be thought of as operations extending $\delta$.
\begin{Def}\label{Def-delta-connection}
Let $\A \xrightarrow{\delta} \Omega$ be a dg derivation and $\E$ an $\A$-module. A $\delta$-connection on $\E$ is a degree $0$, $\k$-linear map of graded $\k$-vector spaces
  \bd
   \nabla :~ \;\E  \rightarrow  \Omega  \otimes_\A \E
  \ed satisfying the following Leibniz rule:
  \be\label{Leib of connections}
   \nabla(ae) = \delta(a) \otimes e + a\nabla(e),\;\;\forall a \in \A , e \in \E.
  \ee
\end{Def}
\begin{Rem}
A connection $\blacktriangledown$ as in Definition~\ref{Def Costello} induces a $\delta$-connection $\nabla$ as in Definition \ref{Def-delta-connection} via the following triangle
  \be\label{delta connection}
   \xymatrix{
   \E  \ar[r]^-{\blacktriangledown} \ar[dr]_-{\nabla } &  \Omega_{\A\mid\k}^1 \otimes_\A \E   \ar[d]^-{\bar{\delta}\otimes_\A \id_{\E}} \\
   &  \Omega \otimes_\A \E .
  }
  \ee
 It follows that $\delta$-connections always exist on projective $\A$-modules.
However, $\delta$-connections do not necessarily arise in this manner.
\end{Rem}

\begin{prop}\label{prop:Atiyah via connection}
Let $\E = (\E, \partial_\A)$ be a dg $\A$-module.\begin{itemize}
\item[1)] For any $\delta$-connection $\nabla $ on $\E$, the degree $1$ element
\bd
   \At_\E^{\nabla}:=[\nabla,\partial_\A] = \nabla  \circ \partial_\A - \partial_\A \circ \nabla  \in \Omega \otimes_\A \End_\A(\E)
\ed
is a cocycle.
\item[ 2)]The cohomology class $[\At^{\nabla }_\E] \in H^1(\A, \Omega \otimes_\A \End_\A(\E))$ coincides with the $\delta$-twisted Atiyah class $\At_\E^\delta$ of $\E$.\end{itemize}
\end{prop}
\bp
The first statement is clear.
It only suffices to prove the second one: Observe that the difference of two $\delta$-connections is a degree zero element in $ \Omega \otimes_\A \End_\A(\E)$. Hence, the cohomology class $[\At_\E^{\nabla }]$ is independent of the choice of $\delta$-connections. We choose a particular $\delta$-connection $\nabla $ induced by a connection $\blacktriangledown$ on $\E$ as in the commutative triangle~\eqref{delta connection}.

Since the map $\bar{\delta} \otimes_\A \id_{\End_\A(\E)}$ defined in Equation~\eqref{Eq:bardeltatensor1}
is a dg $\A$-module morphism, it follows that
 \begin{align*}
   \At_\E^{\nabla } &= [\nabla ,\partial_\A] = [(\bar{\delta} \otimes_\A \id_{\End_\A(\E)}) \circ \blacktriangledown,\partial_\A] = (\bar{\delta} \otimes_\A \id_{\End_\A(\E)})[\blacktriangledown,\partial_\A] = (\bar{\delta} \otimes_\A \id_{\End_\A(\E)})(\At_\E^\blacktriangledown).
 \end{align*}
Passing to the cohomology, we have
 \bd
  [\At_\E^{\nabla }] = (\bar{\delta} \otimes_\A \id_{\End_\A(\E)})(\At_\E) = \At_\E^\delta.
 \ed
\ep

\begin{Def}\label{Def: twisted Atiyah cocycle}
We call $\At_\E^{\nabla }$ the  $\delta$-twisted Atiyah cocycle of $\E$ with respect to the $\delta$-connection $\nabla $ on $\E$.
\end{Def}
Denote the $\A$-dual $\Omega^\vee$ of $\Omega$ by $B$, which is also a dg $\A$-module. Given a $\delta$-connection $\nabla:\;\E  \rightarrow  \Omega \otimes_\A \E$ of an $\A$-module $\E$, the covariant derivation along $b \in \B$ is
$$
\nabla_b:~\E \to \E ,\qquad \nabla_b(e):= \iota_b\nabla(e),\quad~\forall e\in \E.
$$
The $\delta$-twisted Atiyah cocycle $\At_\E^{\nabla }$ could be viewed as a degree $1$ element in $\Hom_\A(\B \otimes_\A \E,\E)$ by setting
\begin{align}\label{Atiyah cocycle}
\At_\E^{\nabla }(b,e) &= (-1)^{\abs{b}}\iota_b\At_\E^{\nabla }(e) = (-1)^{\abs{b}}\iota_b(\nabla (\partial_\A  (e)) - \partial_\A (\nabla (e))) \notag \\
&= -\partial_\A  (\nabla _b e) + \nabla _{\partial_\A (b)}e + (-1)^{\abs{b}}\nabla _b\partial_\A  (e) \notag \\
&= \nabla _{\partial_\A   (b)}e - [\partial_\A  ,\nabla _b](e),
\end{align}
for all $b \in \B $ and $e \in \E$.
Moreover, as $\At_\E^{\nabla }$ is a $1$-cocycle, it is a morphism of dg $\A$-modules, i.e., $\At_\E^{\nabla }\in\Hom^1_{\DG\A}(\B \otimes_\A \E,\E)$.

As an immediate consequence of Proposition~\ref{prop:Atiyah via connection} and Equation~\eqref{Atiyah cocycle}, we have the following
\begin{prop}\label{prop: Atiyah vanish}
Let $\A \xrightarrow{\delta} \Omega$  be a dg derivation and $\E$ a dg $\A$-module. Then the $\delta$-twisted Atiyah class $\At^\delta_\E$ vanishes if and only if there exists a $\delta$-connection $\nabla$ on $\E$ such that the associated twisted Atiyah cocycle $\At_\E^{\nabla }$ vanishes, i.e., the map $\nabla:~\E \rightarrow \Omega \otimes_\A \E$ is compatible with the differentials. In this case,  for all $\partial_\A$-closed elements $b\in \B$ and $ e\in \E$, $\nabla _b e$ is also $\partial_\A$-closed.
\end{prop}

\subsection{Atiyah classes of dg Lie algebroids and Lie pairs}\label{Sec:DGLieAlgebroidsLiePairs}
In this section, we briefly recall Atiyah classes of dg vector bundles with respect to a dg Lie algebroid defined in~\cite{MSX} and Atiyah classes of Lie pairs defined in~\cite{CSX} (see~\cite{CXX} for the equivalence between the two types of Atiyah classes arising from integrable distributions), and show that both of them can be viewed as twisted Atiyah classes.

\subsubsection{Dg Lie algebroids}\label{Ex:DGLiealgebroids}
A dg Lie algebroid can be thought of as a Lie algebroid object in the category of smooth dg manifolds. The precise description is as follows.
\begin{Def}
\label{Def DG Liealgebroid}
  A dg Lie algebroid over a dg manifold $(\M,Q_\M)$ is a quadruple
  $$
  (\L, Q_\L,\rho_\L, [-,-]_{\L}),
  $$
  where
  \begin{itemize}
    \item[1)] $(\L,Q_\L)$ is a dg vector bundle over $(\M,Q_\M)$;
    \item[2)] $(\L,\rho_\L, [-,-]_\L)$ is a graded Lie algebroid over $\M$;
    \item[3)] The anchor map $\rho_\L:~ (\L,Q_\L) \rightarrow (T\M,L_{Q_\M})$ is a morphism of dg vector bundles;
    \item[4)] $Q_\L:~ \Gamma(\L) \rightarrow \Gamma(\L)$ is a derivation with respect to the bracket $[-,-]_\L$, i.e.,
        \bd
        Q_\L([X,Y]_\L) = [Q_\L(X),Y]_\L + (-1)^{\abs{X}}[X,Q_\L(Y)]_\L,\;\;\forall X,Y \in \Gamma(\L).
        \ed
  \end{itemize}
\end{Def}

 Given a dg Lie algebroid $(\L, Q_\L,\rho_\L, [-,-]_{\L})$ and a dg vector bundle $(\mathcal{E},Q_{\mathcal{E}})$ over $(\M,Q_\M)$, Mehta, Sti\'{e}non and Xu constructed the Atiyah class $\At_{\mathcal{E}}$ of $\mathcal{E}$ with respect to $\L$ as follows:
Choose a Lie algebroid $\L$-connection $\nabla^{\mathcal{E}}$ on the vector bundle $\mathcal{E}$, i.e., a degree $0$ $\k$-bilinear map
\bd
 \nabla^{\mathcal{E}}:~ \Gamma(\L) \times \Gamma(\mathcal{E}) \rightarrow \Gamma(\mathcal{E})
\ed
subject to the relations
\begin{align*}
  \nabla^{\mathcal{E}}_{fX}e &= f\nabla^{\mathcal{E}}_{X}e, & \nabla^{\mathcal{E}}_{X}(fe) &= (\rho_\L(X)f)e + (-1)^{\abs{f}\abs{X}} f\nabla^{\mathcal{E}}_{X}e,
\end{align*}
for all $f \in C^\infty(\M), X \in \Gamma(\L)$ and $e \in \Gamma(\mathcal{E})$. There associates a degree $1$ cocycle $\At_{\mathcal{E}}^{\nabla^{\mathcal{E}}} \in \Gamma(\L^\vee \otimes \End(\mathcal{E}))$ defined by
\be\label{Eq: Def of Atiyah class}
 \At_{\mathcal{E}}^{\nabla^{\mathcal{E}}}(X,e) = Q_{\mathcal{E}}(\nabla^{\mathcal{E}}_X e) - \nabla^{\mathcal{E}}_{Q_\L(X)}e - (-1)^{\abs{X}}\nabla^{\mathcal{E}}_X (Q_{\mathcal{E}}e), \;\;\forall X \in \Gamma(\L), e \in \Gamma(\mathcal{E}).
\ee
Its cohomology class $\At_{\mathcal{E}} \in H^1(\Gamma(\L^\vee \otimes \End(\mathcal{E})))$, which is independent of the choice of $\L$-connections, is called the Atiyah class of the dg vector bundle $\mathcal{E}$ with respect to the dg Lie algebroid $\L$~\cite{MSX}.

Meanwhile, there associates a $(\Gamma(\L^\vee),Q_{\L^\vee})$-valued derivation of the cdga $(C^\infty(\M),Q_\M)$ defined by
  \begin{equation}\label{Eq:deltaL}
   \delta_\L: \; C^\infty(\M)  \xrightarrow{ d_{dR}} \Omega^1(\M)\xrightarrow{\rho_\L^\vee} \Gamma(\L^\vee),
  \end{equation}
where $Q_{\L^\vee}$ is induced from the differential $Q_\L$ on $\L$. The fact that $\delta_\L$ commutes with the two differentials $Q_{\M}$ and $Q_{\L^\vee}$ follows from (3) of Definition \ref{Def DG Liealgebroid}.
The section space of a dg vector bundle $\mathcal{E}$ gives rise to a dg $(C^\infty(\M),Q_\M)$-module $\E:=(\Gamma(\mathcal{E}),Q_\mathcal{E})$.
It is obvious that a $\delta_\L$-connection $\nabla^{\delta_\L}$ on $\E= \Gamma(\mathcal{E})$ is equivalently to a Lie algebroid $\L$-connection $\nabla^\L$ on the graded vector bundle $\mathcal{E}$.
Comparing Equations~\eqref{Atiyah cocycle} and~\eqref{Eq: Def of Atiyah class}, we have the following
\begin{prop}\label{prop: DG Lie algebroid}
  The Atiyah class $\At_{\mathcal{E}}$ of the dg vector bundle $\mathcal{E}$ with respect to the dg Lie algebroid $\L$  coincides, up to a minus sign, with the twisted Atiyah class $\At_\E^{\delta_\L}$ of the dg $(C^\infty(\M),Q_\M)$-module $\E=(\Gamma(\mathcal{E}),Q_\mathcal{E})$, where  the dg derivation $\delta_\L$ is given by Equation \eqref{Eq:deltaL}.
\end{prop}

\subsubsection{Lie pairs}\label{Ex: Lie pair}
By a Lie pair $(L,A)$, we mean two Lie algebroids $L$ and $A$ over the same smooth manifold $M$ such that $A \subset L$ is a Lie subalgebroid.  The quotient bundle $B=L/A$ carries a natural flat (Lie algebroid) $A$-connection, called the Bott $A$-module structure.

Let us recall the Atiyah class of the Lie pair $(L,A)$ defined in~\cite{CSX}.
First of all, there is a short exact sequence of vector bundles over $M$,
\be\label{SES}
 0 \rightarrow A \xrightarrow{~i~} L \xrightarrow{\pr_B} B \rightarrow 0.
\ee
Choose a splitting of Sequence~\eqref{SES}, i.e., a vector bundle injection  $j:~ B \rightarrow L$, which determines a  bundle projection $\pr_A:~ L \rightarrow A$ such that
\begin{align*}
\pr_A \circ i &= \id_A, & \pr_B \circ j &= \id_B, & i \circ \pr_A + j \circ \pr_B &= \id_L.
\end{align*}
Using this splitting, one could identify $L$ with $  A\oplus B$.
Meanwhile, for each $A$-module $(E,\partial_A^E)$, where $E$ is a vector bundle over $M$ and $\partial_A^E$ is a flat $A$-connection on $E$, choose  an $L$-connection $\nabla^L$ on $E$ extending the given flat $A$-connection.
Then there associates a $1$-cocycle $\alpha_E^{\nabla^L} \in \Gamma(A^\vee \otimes B^\vee \otimes \End(E))$, called the Atiyah cocycle, of the Lie algebroid $A$ valued in the $A$-module $B^\vee \otimes \End(E)$:
\begin{align*}
 \alpha_E^{\nabla^L}(a,b)e &:= \nabla_{a}\nabla^L_{j(b)}e - \nabla^L_{j(b)}\nabla_{a}e - \nabla^L_{[a,j(b)]}e,
\end{align*}
for all $a \in \Gamma(A), b \in \Gamma(B)$ and $e \in \Gamma(E)$. 
The cohomology class
\bd
\alpha_{E} = [\alpha_E^{\nabla^L}] \in H_{\CE}^1(A,B^\vee \otimes \End(E))
\ed
does not depend on the choice of $j$ and $\nabla^L$, and is called the Atiyah class of the $A$-module $E$ with respect to the Lie pair $(L,A)$.

From the Lie pair $(L,A)$, we get a cdga $\Omega^\bullet_A = (\Gamma(\wedge^\bullet A^\vee),d_{A})$, and a dg $\Omega^\bullet_A$-module
$\Omega_A^\bullet(B^\vee):= (\Gamma(\wedge^\bullet A^\vee \otimes B^\vee),\partial_A )$, where $\partial_A$ is the $A$-module structure dual to the Bott $A$-module structure on $B$.
Here the degree convention is that $\Gamma({B^\vee})$ concentrates in degree zero.

Fixing a splitting $j$ of Sequence~\eqref{SES}, we construct an
 $\Omega_A^\bullet(B^\vee)$-valued derivation $\delta_j$ of $\Omega^\bullet_A$, i.e., a map
\begin{equation}\label{Eqt:deltaofLA}
\delta_j:~\Gamma(\wedge^\bullet A^\vee ) \to \Gamma(\wedge^\bullet A^\vee \otimes B^\vee).
\end{equation}
As a degree zero derivation of the graded $\k$-algebra $\Gamma(\wedge^\bullet A^\vee)$, $\delta_j$ is fully determined by its action on its generators, i.e. elements in $C^\infty(M)$ and $\Gamma(A^\vee)$ --- Define
 \begin{align*}
    \delta_j:~ \;&C^\infty(M) \xrightarrow{d_L}   \Gamma(L^\vee) \xrightarrow{j^\vee} \Gamma(B^\vee), \\
    \delta_j:~ \;&\Gamma( A^\vee) \xrightarrow{\pr_A^\vee} \Gamma( L^\vee) \xrightarrow{d_L} \Gamma(\wedge^2 L^\vee) \xrightarrow{ } \Gamma( L^\vee \otimes L^\vee)  \xrightarrow{i^\vee \otimes j^\vee} \Gamma( A^\vee \otimes B^\vee),
 \end{align*}
where $d_L:~ \Gamma(\wedge^\bullet L^\vee) \rightarrow \Gamma(\wedge^{\bullet+1}L^\vee)$ is the Chevalley-Eilenberg differential of the Lie algebroid $L$. A straightforward verification shows that $\delta_j$ is compatible with the differentials and thus is an $\Omega^\bullet_A(B^\vee)$-valued dg derivation of $\Omega^\bullet_A$.

Note that $\delta_j$ depends on a choice of a splitting $j$ of Sequence~\eqref{SES}. However, we have
\begin{prop}\label{prop: splitting}
 The $\Omega_A^\bullet(B^\vee)$-valued dg derivations $\delta_j$ of $\Omega_A^\bullet$ associated with different splittings of Sequence~\eqref{SES} are homotopic to each other.
\end{prop}
\bp
 Given two splittings $j$ and $j^\prime$ of Sequence~\eqref{SES}, their difference is a bundle map $j^\prime-j:~ B \rightarrow A$. Define a degree $(-1)$ derivation $h:~\Gamma(\wedge^{\bullet}A^\vee) \rightarrow \Gamma(\wedge^{\bullet-1}A^\vee \otimes B^\vee)$ by setting
 $$
  h|_{C^\infty(M)}=0,\qquad h|_{\Gamma(A^\vee)}=(j^\prime-j)^\vee.
 $$
 It follows from direct verifications that
 \bd
  \delta_{j^\prime} - \delta_j = [\partial_A,h]: \Gamma(\wedge^{\bullet}A^\vee) \rightarrow \Gamma(\wedge^{\bullet}A^\vee \otimes B^\vee). \ed
This proves that $\delta_j\sim \delta_{j^\prime}$.
\ep

Let $(E,\partial_A^E)$ be an $A$-module. There induces a dg $\Omega^\bullet_A$-module  $\E := (\Gamma(\wedge^\bullet A^\vee \otimes E),\partial_A^E)$.
\begin{prop}\label{prop:Lie pairs}
  The Atiyah class $\alpha_E$ of the $A$-module $E$ with respect to the Lie pair $(L,A)$ coincides with the twisted  Atiyah class $\At^{\delta_j}_\E$ of the dg $\Omega_A^\bullet$-module $\E$, where the dg derivation $\delta_j$ is given as in Equation \eqref{Eqt:deltaofLA}.
\end{prop}
\bp
First of all, the spaces where the two Atiyah classes live are exactly the same, i.e.,
\bd
   (\alpha_E \in)~ H^1_{\CE}(A,B^\vee \otimes \End(E)) = H^1(\Omega_A^\bullet, \Omega_A^\bullet(B^\vee) \otimes_{\Omega_A^\bullet} \End_{\Omega_A^\bullet}(\E)) ~(\ni \At^{\delta_j}_\E).
\ed
According to Proposition \ref{prop:Atiyah via connection}, to find the twisted Atiyah class $\At^{\delta_j}_\E$, one may use   a $\delta_j$-connection $\nabla^{\delta_j}$ on $\E$, which is determined by its restriction on $\Gamma(E)$:
\bd
 \nabla^{\delta_j}\!\mid_E:~ \Gamma(E) \rightarrow \Gamma(B^\vee) \otimes \Gamma(E).
\ed
This is equivalent to an $L$-connection $\nabla^L$ on $E$ extending the given flat $A$-connection by setting
\be\label{Eq:Landdeltaconnection}
 \nabla^L_{a+b} = \nabla_a + \nabla^{\delta_j}_b\!\mid_E,\qquad \forall a+b\in L\cong A\oplus B.
\ee
The two associated Atiyah cocycles coincide by straightforward computations, i.e., $\At_\E^{\nabla^{\delta_j}} = \alpha_E^{\nabla^L}$.
\ep
As a consequence of Propositions~\ref{prop: DG Lie algebroid} and~\ref{prop:Lie pairs}, both Atiyah classes of dg Lie algebroids and those of Lie pairs arise from Atiyah classes of cdgas.
In particular, we have
\begin{Cor}
  Let $A$ be a Lie algebroid and $E$ an $A$-module. Denote by $\mathcal{E}$ the corresponding dg vector bundle over $(A[1],d_A)$. 
  If the Atiyah class of the dg vector bundle $\mathcal{E}$ with respect to the dg Lie algebroid $T(A[1])$ vanishes, then the Atiyah class of $E$ with respect to any Lie pair $(L,A)$ vanishes.
\end{Cor}

\subsection{Functoriality}
We now study functorial properties of twisted Atiyah classes. 
Let $\CAThtyAmod$ denote the homology category of dg $\A$-modules: Objects in $\CAThtyAmod$ are dg $\A$-modules, and morphisms in $\CAThtyAmod$ are dg $\A$-module morphisms modulo homotopy~\cite{Keller}.

Let $\A \xrightarrow{\delta} \Omega$ be a dg derivation.
For each object $\E$ in $\CAThtyAmod$, by Definition~\ref{Def: twisted Atiyah class}, the twisted Atiyah class
$$
\At^\delta_\E \in H^1(\A, \Omega \otimes_\A \End_\A(\E))\cong
\Hom^{1}_{\CAThtyAmod}(\E,\Omega \otimes_\A \E)
$$
is a degree $1$ morphism in the category $\CAThtyAmod$. This identification defines a functorial transformation.
In fact, when the dg derivation $\delta$ is fixed, the $\delta$-twisted Atiyah class is a functorial transformation on $\CAThtyAmod$ from the identity functor $\id$ to the tensor functor $\Omega \otimes_\A -$:
\begin{prop}\label{Prop: functorial 1}
Let $\E$ and $\F$ be dg $\A$-modules, $\lambda \in \Hom_{\CAThtyAmod}(\E,\F)$. The following diagram commutes in the category $\CAThtyAmod$:
\begin{equation}\label{Eqt:twisteddiagram}
\xymatrix{
   \E \ar[r]^-{\At^\delta_\E} \ar[d]_-{(-1)^{\abs{\lambda}}\lambda} & \Omega \otimes_\A \E \ar[d]^-{\id_\Omega \otimes_\A \lambda} \\
   \F \ar[r]^-{\At^\delta_\F} & \Omega \otimes_\A \F.
}
\end{equation}
\end{prop}
\bp
Let us first show the non-twisted case. Namely, the following diagram commutes in $\CAThtyAmod$:
\begin{equation}\label{Eqt:nontwisteddigram}
\xymatrix{
   \E \ar[r]^-{\At_\E} \ar[d]_-{(-1)^{\abs{\lambda}}\lambda} & \Omega_{\A\mid\k}^1 \otimes_\A \E \ar[d]^-{\id  \otimes_\A \lambda} \\
   \F \ar[r]^-{\At_\F} & \Omega_{\A\mid\k}^1 \otimes_\A \F.
}
\end{equation}
In fact, this can be directly verified. We choose connections $\blacktriangledown^\E$ and $\blacktriangledown^\F$, respectively, on $\E$ and $\F$. For simplicity, they are both denoted by $\blacktriangledown$. Then
\begin{align*}
  &\quad(\id  \otimes_\A \lambda) \circ \At^{\blacktriangledown}_\E - (-1)^{\abs{\lambda}}\At_\F^{\blacktriangledown} \circ \lambda\\
  &= (\id_\Omega \otimes_\A \lambda) \circ (\blacktriangledown \circ \partial_\A   - \partial_\A  \circ \blacktriangledown) - (-1)^{\abs{\lambda}}(\blacktriangledown \circ \partial_\A - \partial_\A  \circ \blacktriangledown) \circ \lambda\\
  &=  ((\id  \otimes_\A \lambda) \circ \blacktriangledown - \blacktriangledown \circ \lambda) \circ \partial_\A   - (-1)^{\abs{\lambda}}\partial_\A  \circ ((\id  \otimes_\A \lambda) \circ \blacktriangledown - \blacktriangledown \circ \lambda).
\end{align*}
The map $(\id  \otimes_\A \lambda) \circ \blacktriangledown - \blacktriangledown \circ \lambda:~
\E  \to  \Omega_{\A\mid\k}^1  \otimes_{\A } \F  $ is actually $  {\A }$-linear, by direct verifications. Thus the two maps $(\id  \otimes_\A \lambda) \circ \At^{\blacktriangledown}_\E$  and $ (-1)^{\abs{\lambda}}\At_\F^{\blacktriangledown} \circ \lambda$ are only differed by an exact term. Composing with the dg $\A$-module morphism $\bar{\delta}:~\Omega_{\A\mid\k}^1 \to \Omega$ induced from the dg derivation $\delta$, we accomplish a commutative diagram in $\CAThtyAmod$:
$$
\xymatrix{
   \E \ar[r]^-{\At_\E} \ar[d]_-{(-1)^{\abs{\lambda}}\lambda} &  \Omega_{\A\mid\k}^1 \otimes_\A \E \ar[d]^-{\id  \otimes_\A \lambda} \ar[r]^-{\bar{\delta}\otimes_\A \id} & \Omega  \otimes_\A \E \ar[d]^-{\id  \otimes_\A \lambda}\\
   \F \ar[r]^-{\At_\F} & \Omega_{\A\mid\k}^1 \otimes_\A \F  \ar[r]^-{\bar{\delta}\otimes_\A \id} &\Omega  \otimes_\A \F.
}
$$
\ep

Now we study how Atiyah classes vary when twisted by different dg derivations. So we need the category of dg derivations, denoted by $\CATderivationA$, whose objects are dg derivations $\A \xrightarrow{\delta} \Omega$ as in Definition \ref{Def:DGderivations}, and whose morphisms are defined as follows:
\begin{Def}\label{Def:morphismDerA}
  A morphism $\phi$ from $\A \xrightarrow{\delta} \Omega$ to $\A \xrightarrow{\delta^\prime} \Omega^\prime$ is a morphism $\phi: \Omega \rightarrow \Omega^\prime$ of dg $\A$-modules such that
  \bd
    \delta^\prime = \phi \circ \delta: \A \rightarrow \Omega^\prime.
  \ed
\end{Def}

Now let us fix a dg $\A$-module $\E$. We have the constant functor $(- \mapsto \E)$ and the tensor functor $- \otimes_\A \E$, both from the category $\CATderivationA$ of dg derivations to the homology category $\CAThtyAmod$ of dg $\A$-modules.
The Atiyah class is a functorial transformation from $(-\mapsto \E)$ to $- \otimes_\A \E$:
\begin{prop}\label{prop:functor2}
  Given a morphism $\phi:~ (\A \xrightarrow{\delta} \Omega) \rightarrow (\A \xrightarrow{\widetilde{\delta}} \widetilde{\Omega})$ of dg derivations  and a dg $\A$-module $\E$, let $\At^\delta_\E$ and $\At^{\tilde{\delta}}_\E$ be the Atiyah classes of $\E$ twisted, respectively, by $\A \xrightarrow{\delta} \Omega$ and $\A \xrightarrow{\widetilde{\delta}} \widetilde{\Omega}$. Then the following diagram commutes in $\CAThtyAmod$:
  \bd
  \xymatrix{
     \E \ar[r]^-{\At^\delta_\E} \ar[d]_-{\id_\E} & \Omega \otimes_\A \E \ar[d]^-{\phi \otimes_\A \id_\E} \\
     \E \ar[r]^-{\At^{\tilde{\delta}}_\E} & \widetilde{\Omega} \otimes_\A \E.
  }
  \ed
\end{prop}
The proof is easy and thus omitted. Combining the previous two propositions, we have
\begin{Thm}\label{Thm: Functoriality}
  With the same assumptions as in Propositions~\ref{Prop: functorial 1} and~\ref{prop:functor2}, the following diagram in $\CAThtyAmod$ commutes:
  \bd
\xymatrix{
   \E \ar[r]^-{\At^{\delta}_\E} \ar[d]_-{(-1)^{\abs{\lambda}}\lambda} & \Omega \otimes_\A \E \ar[d]^-{\phi \otimes_\A \lambda} \\
   \F \ar[r]^-{\At^{\tilde{\delta}}_\F} & \widetilde{\Omega} \otimes_\A \F.
}
\ed
\end{Thm}

\section{The Kapranov functor}\label{Sec:MainResults}
In this section, we explore higher algebraic structures, called Kapranov Leibniz$_\infty[1]$ algebras, induced from a dg derivation of a cdga $\A$. Our main goal is to show that there exists a contravariant functor, called the Kapranov functor, from the category of dg derivations to the category of Leibniz$_\infty[1]$-algebras over $\A$.

\subsection{Leibniz$_\infty[1]$ algebras}
We recall some basic notions of homotopy Leibniz algebras (c.f.\cites{AP,CSX}). In what follows, all tensor products $\otimes$ without adoration are assumed to be over $\k$.
\begin{Def}\label{Def:Leibnizinfinity}
  A Leibniz$_\infty[1]$ algebra (over $\k$) is a graded $\k$-vector space $V = \oplus_{n \in \Z}V^n$,  together with a sequence $\{\lambda_k:~ \otimes^k V \rightarrow V\}_{k \geq 1}$ of degree $1$, $\k$-multilinear maps satisfying
  \begin{align}\label{Leib}
    &\sum_{i+j=n+1}\sum_{k=j}^{n}\sum_{\sigma \in \sh(k-j,j-1)}\epsilon(\sigma)(-1)^{\abs{v_{\sigma(1)}}+\cdots+ \abs{v_{\sigma(k-j)}}} \notag \\
    &\lambda_{i}(v_{\sigma(1)},\cdots,v_{\sigma(k-j)},\lambda_{j}(v_{\sigma(k-j+1)}, \cdots,v_{\sigma(k-1)}, v_k), v_{k+1},\cdots, v_n)=0,
  \end{align}
  for all $n \geq 1$ and all homogeneous elements $v_i \in V$, where $\sh(p,q)$ denotes the set of $(p,q)$-shuffles ($p,q \geq 0$), and $\epsilon(\sigma)$ is the Koszul sign of  $\sigma$.
\end{Def}

\begin{Def}
A  morphism of Leibniz$_\infty[1]$ algebras from $(V,\{\lambda_k\}_{k\geq1})$ to $(V',\{\lambda'_k\}_{k\geq1})$ is a sequence $\{f_k:~ V^{\otimes k} \rightarrow V'\}_{k \geq 1}$ of degree $0$, $\k$-multilinear maps, satisfying the following compatibility condition:
\begin{align}\label{morphism equation}
  &\quad \sum_{k+p \leq n-1}\sum_{\sigma \in \sh(k,p)}\epsilon(\sigma)(-1)^{\dagger^{\sigma}_k}f_{n-p}(b_{\sigma(1)},\cdots,b_{\sigma(k)}, \lambda_{p+1}(b_{\sigma(k+1)}, \cdots,b_{\sigma(k+p)},b_{k+p+1}),\cdots,b_{n}) \notag \\
  &= \sum_{q \geq 1}\sum_{\substack{I^1\cup\cdots\cup I^q = \mathbb{N}^{(n)} \\ I_1,\cdots,I_q \neq \emptyset \\ i_{\abs{I^1}}^1 < \cdots < i^q_{\abs{I^q}}}}\epsilon(I^1,\cdots,I^q) \lambda_q'(f_{\abs{I^1}}(b_{I^1}), \cdots,f_{\abs{I^q}}(b_{I^q})),
 \end{align}
for all $n \geq 1$, where $\dagger^{\sigma}_k = \sum_{i=1}^k\abs{b_{\sigma(i)}}$, $I^j = \{i^j_1 < \cdots < i^j_{\abs{I^j}}\} \subset \mathbb{N}^{(n)} = \{1,\cdots,n\}$, and $(b_{I^j}) = (b_{i^j_1},\cdots,b_{i^j_{\abs{I^j}}})$ for all $1\leq j \leq q$.
\end{Def}
In the above definition, the first component $f_1:~ (V,\lambda_1) \rightarrow (V',\lambda_1')$, called the tangent morphism, is a morphism of cochain complexes. We call the Leibniz$_\infty[1]$ morphism $f_\bullet:~ (V,\lambda_\bullet) \rightarrow (V',\lambda'_\bullet)$ a quasi-isomorphism (resp. an isomorphism) if $f_1$ is a quasi-isomorphism (resp. an isomorphism). In fact, there is a standard way to find its quasi-inverse (resp. inverse)
$f^{-1}_\bullet:~ (V',\lambda'_\bullet) \rightarrow (V,\lambda_\bullet)$ (see \cite{AP}).

\begin{Def}\label{Def:Leibnizinfinitymodule}Let $(V,\{\lambda_k\}_{k\geq1})$ be a Leibniz$_\infty[1]$ algebra.
  A $(V,\{\lambda_k\}_{k\geq1})$-module  is a graded $\k$-vector space $W$ together with a sequence $\{\mu_k:~ V^{\otimes(k-1)} \otimes W \rightarrow W\}_{k \geq 1}$ of degree $1$, $\k$-multilinear maps satisfying the identities
  \begin{align*}
    &\sum_{i+j=n+1}\sum_{k=j}^{n}\sum_{\sigma \in \sh(k-j,j-1)}\epsilon(\sigma)(-1)^{\dagger^\sigma_{k-j}}  \notag\\
    &\qquad\mu_{i}(v_{\sigma(1)},\cdots,v_{\sigma(k-j)},\lambda_{j}(v_{\sigma(k-j+1)}, \cdots,v_{\sigma(k-1)}, v_k), v_{k+1},\cdots, v_{n-1},w) \notag \\
    &\qquad+\sum_{1\leq j \leq n}\sum_{\sigma \in \sh(k,j)}\epsilon(\sigma)(-1)^{\dagger^\sigma_{n-j}} \mu_{i}(v_{\sigma(1)},\cdots,v_{\sigma(n-j)},\mu_{j}(v_{\sigma(n-j+1)}, \cdots,v_{\sigma(n-1)}, w))=0,
  \end{align*}
  for all $n \geq 1$ and all homogeneous vectors $v_1,\cdots,v_{n-1} \in V, w \in W$, where $\dagger^\sigma_j = \abs{v_{\sigma(1)}}+\cdots+ \abs{v_{\sigma(j)}}$ for all $j \geq 0$.
 \end{Def}

\begin{Def}
Let $(V,\{\lambda_k\}_{k\geq1})$ be a Leibniz$_\infty[1]$ algebra.
   A morphism of $(V,\{\lambda_k\}_{k\geq1})$-modules from $(W,\{\mu_k\}_{k\geq1})$ to $(W',\{\mu'_k\}_{k\geq1})$  is a sequence $\{\psi_k:~ V^{\otimes(k-1)} \otimes W \rightarrow W'\}_{k \geq 1}$ of degree $0$,  $\k$-multilinear maps satisfying the identity
   \begin{align*}\label{Eq:Leib module morphism}
    &\sum_{i+j=n+1}\sum_{k=j}^{n}\sum_{\sigma \in \sh(k-j,j-1)}\epsilon(\sigma)(-1)^{\abs{v_{\sigma(1)}}+\cdots+ \abs{v_{\sigma(k-j)}}}  \notag\\
    &\quad\psi_{i}(v_{\sigma(1)},\cdots,v_{\sigma(k-j)},\lambda_{j}(v_{\sigma(k-j+1)}, \cdots,v_{\sigma(k-1)}, v_k), v_{k+1},\cdots, v_{n-1},w) \notag \\
    &\qquad+\sum_{1\leq j \leq n}\sum_{\sigma \in \sh(k,j)}\epsilon(\sigma)(-1)^{\abs{v_{\sigma(1)}}+\cdots+ \abs{v_{\sigma(n-j)}}} \psi_{i}(v_{\sigma(1)},\cdots,v_{\sigma(n-j)},\mu_{j}(v_{\sigma(n-j+1)},\cdots, v_{\sigma(n-1)}, w)) \notag\\
   &= \sum_{p \geq 0}\sum_{\substack{I^1\cup\cdots\cup I^{p+1} = \mathbb{N}^{(n-1)} \\ I^1,\cdots,I^{p+1} \neq \emptyset \\ i_{\abs{I^1}}^1 < \cdots < i^{p+1}_{\abs{I^{p+1}}}}}\epsilon(I^1,\cdots,I^{p+1}) \mu'_{p+1}(\lambda_{\abs{I^1}}(v_{I^1}),\cdots,\lambda_{\abs{I^p}}(v_{I^p}), \psi_{\abs{I^{p+1}}+1}(v_{I^{p+1}},w)),
  \end{align*}
  for each $n \geq 1$ and all homogeneous vectors $v_1,\cdots,v_{n-1} \in V, w \in W$. Here $I^j = \{i^j_1 < \cdots < i^j_{\abs{I^j}}\} \subset \mathbb{N}^{(n-1)} = \{1,\cdots,n-1\}$, and $(v_{I^j}) = (v_{i^j_1},\cdots,v_{i^j_{\abs{I^j}}})$ for all $1 \leq j \leq p+1$.
\end{Def}

In this note, we are particularly interested in Leibniz$_\infty[1]$ algebras over a cdga $\A$ (or Leibniz$_\infty[1]$ $\A$-algebras).
\begin{Def}\label{Def:LeibnizinfinityA}
A Leibniz$_\infty[1]$ $\A$-algebra is a Leibniz$_\infty[1]$ algebra $(V,\{\lambda_k\}_{k\geq1})$ (in the category of Leibniz$_\infty[1]$ algebras over $\k$) such that the cochain complex $(V,\lambda_1)$ is a dg $\A$-module and all higher brackets $\lambda_k: \otimes^k V \rightarrow V$ ($k\geq 2$) are  $\A$-multilinear.

A morphism of Leibniz$_\infty[1]$  $\A$-algebras from $(V,\{\lambda_k\}_{k\geq1})$ to $(V',\{\lambda'_k\}_{k\geq1})$ is a
morphism $\{f_k:V^{\otimes k}\rightarrow V'\}_{k \geq 1}$ (in the category of Leibniz$_\infty[1]$ algebras over $\k$) such that all structure maps $\{f_k\}_{k\geq 1}$ are  $\A$-multilinear. In particular, its tangent morphism $f_1:~ (V,\lambda_1) \rightarrow (V^\prime,\lambda^\prime_1)$ is a dg $\A$-module morphism.

Such a morphism $f_\bullet:~ (V,\lambda_\bullet) \rightarrow (V',\lambda'_\bullet)$ is called a quasi-isomorphism (resp. an isomorphism) if its tangent morphism $f_1$ is a quasi-isomorphism (resp. an isomorphism) of dg $\A$-modules.
\end{Def}
Denote the category of Leibniz$_\infty[1]$ $\A$-algebras by $\CATLeibnizoneA$. It is a subcategory of the category of Leibniz$_\infty[1]$ algebras over $\k$.
There are analogous notions of modules of a Leibniz$_\infty[1]$ $\A$-algebra and their morphisms:
\begin{Def}\label{Def:LeibnizinfinityAmodule}
Let $(V,\{\lambda_k\}_{k\geq1})$ be a Leibniz$_\infty[1]$ $\A$-algebra.
A $(V,\{\lambda_k\}_{k\geq1})$ $\A$-module is a $(V,\{\lambda_k\}_{k\geq1})$-module $(W,\{\mu_k\}_{k\geq1})$ (in the category of Leibniz$_\infty[1]$ modules over $\k$)   such that $(W,\mu_1)$ is a DG $\A$-module and all higher structure maps $\{\mu_k\}_{k\geq 2}$ are  $\A$-multilinear.

A morphism of $(V,\{\lambda_k\}_{k\geq1})$ $\A$-modules from $(W,\{\mu_k\}_{k\geq1})$ to $(W',\{\mu'_k\}_{k\geq1})$ is a morphism of $(V,\{\lambda_k\}_{k\geq1})$-modules $\{\psi_k:~ V^{\otimes(k-1)} \otimes W \rightarrow W'\}_{k \geq 1}$ (in the category of  Leibniz$_\infty[1]$ modules over $\k$) such that all maps $\{\psi_k\}_{k\geq 1}$ are  $\A$-multilinear.
\end{Def}
It follows that the collection of $(V,\{\lambda_k\}_{k\geq1})$ $\A$-modules and their morphisms form a category.

\subsection{The Kapranov functor}
In this section, we generalize Kapranov's construction of an $L_\infty$ algebra structure~\cite{Kap} and Chen-Sti\'{e}non-Xu's construction of a Leibniz$_\infty[1]$ algebra structure~\cite{CSX} in the setting of a dg derivation $\A \xrightarrow{\delta} \Omega$ of a cdga $\A$.

\subsubsection{Kapranov Leibniz$_\infty[1]$ algebras}
Let $\E$ be a graded $\A$-module with a $\delta$-connection $\nabla$. For each homogeneous $b \in \B$, there is a degree $\abs{b}$ derivation on the reduced tensor algebra $T(\E)$ (over $\A$) defined by
\bd
 \nabla_b(e_1 \otimes \cdots \otimes e_n) = \sum_{i=1}^{n}(-1)^{\abs{b}\ast_{i-1}} e_1 \otimes \cdots \nabla_b e_i \otimes \cdots e_n,
\ed
for all homogeneous $e_i \in \E$, where $\ast_{i} = \sum_{j=1}^i\abs{e_j}$.

Let $\E$ and $\F$ be two graded $\A$-modules with $\delta$-connections $\nabla^{\E}$ and $\nabla^{\F}$, respectively.
For $b \in \B$ and $\lambda \in \Hom_{\A }(\E,\F)$, there associates the derivation
\bd
\nabla_b(\lambda) = [\nabla_b,\lambda] = \nabla^{\F }_b \circ \lambda - (-1)^{\abs{b}\abs{\lambda}}\lambda \circ \nabla^{\E }_b:~ \E  \rightarrow \F .
\ed
It follows from a direct verification that $\nabla_b(\lambda) \in \Hom_{\A}(\E ,\F)$.

Choose a $\delta$-connection on $\B$. There associates a sequence of degree $1$ maps
 $ \mathcal{R}^\nabla_k: \B^{\otimes k} \rightarrow \B, k \geq 1 $
defined as follows:
\begin{itemize}
  \item $\mathcal{R}^\nabla_1 = \partial_\A:~ \B \rightarrow \B$;
  \item $\mathcal{R}^\nabla_2$ is specified by the associated twisted Atiyah cocycles $\At_\B^\nabla$; 
  \item $\{\mathcal{R}^\nabla_{k+1}:~ \B^{\otimes(k+1)} \rightarrow \B\}_{k\geq2}$ are defined recursively by 
      $\mathcal{R}^\nabla_{k+1}=\nabla(\mathcal{R}^\nabla_{k})$. 
     Explicitly, we have
  \begin{align}\label{Rnabla}
   &\mathcal{R}^\nabla_{k+1}(b_0,b_1,\cdots,b_{k}) = (-1)^{\abs{b_0}}[\nabla_{b_0},\mathcal{R}^\nabla_k] (b_1,\cdots,b_k), \;\;\forall b_i \in \B.
  \end{align}
\end{itemize}

\begin{prop}\label{Thm: leibniz infty}
The $\A$-module $\B$, together with the sequence of operators $\{\mathcal{R}^\nabla_k\}_{k\geq 1}$, is a Leibniz$_\infty[1]$ $\A$-algebra.
\end{prop}
\bp
The $2$-bracket $\mathcal{R}^\nabla_2$ is the twisted Atiyah cocycles $\At_\B^\nabla$, which is certainly $\A$-bilinear. By the recursive construction of higher brackets $\mathcal{R}^\nabla_{k+1}$ ($k\geq 2$) in Equation \eqref{Rnabla}, they are all $\A$-multilinear as well.
So it suffices to verify that $\{\mathcal{R}^\nabla_k\}_{k\geq1}$ satisfies Equation~\eqref{Leib}. We argue by induction.

The $n=1$ case follows from the fact that $\mathcal{R}^\nabla_1 = \partial_\A$ is a differential, and the $n=2$ case follows from the fact that the $\delta$-twisted Atiyah cocycle $\At^\nabla_\B$ is a $\partial_\A$-cocycle.

Now assume that the identity~\eqref{Leib} holds for some $n \geq 2$, i.e.,
  \begin{align}\label{casen}
    &\quad -\partial_\A(\mathcal{R}^\nabla_n)(b_1,\cdots,b_n) = -\partial_\A(\mathcal{R}^\nabla_n(b_1,\cdots,b_n))  + \sum_{i=1}^n(-1)^{\ast_{i-1}} \mathcal{R}^\nabla_n(b_1,\cdots,\partial_\A b_i,\cdots,b_{n}) \notag \\
    &=\sum_{i,j \geq 2, i+j=n+1}\sum_{k=j}^{n}\sum_{\sigma \in \sh(k-j,j-1)}\epsilon(\sigma)(-1)^{\abs{b_{\sigma(1)}} + \cdots + \abs{b_{\sigma(k-j)}}} \notag \\
    &\qquad\qquad\mathcal{R}^\nabla_{i}(b_{\sigma(1)},\cdots,b_{\sigma(k-j)}, \mathcal{R}^\nabla_{j} (b_{\sigma(k-j+1)},\cdots,b_{\sigma(k-1)}, b_k), b_{k+1},\cdots, b_n).
  \end{align}
Consider the $(n+1)$ case: We first compute
\begin{align}\label{eq1}
    -\partial_\A(\mathcal{R}^\nabla_{n+1})(b_0,b_1,\cdots,b_n) &= \partial_\A ([\mathcal{R}^\nabla_{n},\nabla_{b_0}])(b_1,\cdots,b_n) + [\mathcal{R}^\nabla_n,\nabla_{\partial_\A b_0}](b_1,\cdots,b_n) \notag \\
    &= [\partial_\A(\mathcal{R}^\nabla_n),\nabla_{b_0}](b_1,\cdots,b_n) + [\mathcal{R}^\nabla_n,\mathcal{R}^\nabla_2(b_0,-)](b_1,\cdots,b_n).
\end{align}
Here we have used the recursive definition~\eqref{Rnabla} in the first equality and Equation~\eqref{Atiyah cocycle} in the second one.

We introduce
\bd
 b^{(i)}_k = \begin{cases}
   b_k, &\;\;\;\text{if}\; k \neq i \\
 \nabla_{b_0}b_i, &\;\;\;\text{if}\; k = i.
 \end{cases}
 \ed
Then the first summand in Equation~\eqref{eq1} is, 
\begin{align}\label{eq2}%
  &\quad [\partial_\A(\mathcal{R}^\nabla_n),\nabla_{b_0}](b_1,b_2,\cdots,b_{n}) \notag\\
  &=\sum_{i=1}^{n}(-1)^{\abs{b_0}\ast_{i-1}}\partial_\A(\mathcal{R}^\nabla_n) (b_1,\cdots,\nabla^\B_{b_0}b_i, \cdots,b_{n}) - \nabla_{b_0}(\partial_\A(\mathcal{R}^\nabla_n)(b_1,\cdots,b_{n})) \notag \;\\&\qquad\qquad \text{by assumption~\eqref{casen}} \notag \\
  &= -\sum_{i=1}^n\sum_{p,q \geq 2, p+q=n+1}\sum_{k=q}^{n}\sum_{\sigma \in \sh(k-q,q-1)}\epsilon(\sigma)(-1)^{\abs{b_0}\ast_{i-1}+\abs{b^{(i)}_{\sigma(1)}} + \cdots+\abs{b^{(i)}_{\sigma(k-q)}}} \notag \\
    &\qquad\qquad\mathcal{R}^\nabla_{p}(b^{(i)}_{\sigma(1)},\cdots, b^{(i)}_{\sigma(k-q)}, \mathcal{R}^\nabla_{j}(b^{(i)}_{\sigma(k-q+1)}, \cdots,b^{(i)}_{\sigma(k-1)}, b^{(i)}_k), b^{(i)}_{k+1},\cdots, b^{(i)}_n) \notag\\
  &\quad +\sum_{i,j \geq 2, i+j=n+1}\sum_{k=j}^{n}\sum_{\sigma \in \sh(k-j,j-1)}\epsilon(\sigma)(-1)^{\abs{b_{\sigma(1)}}+\cdots+\abs{b_{\sigma(k-j)}}} \notag \\
    &\qquad\qquad\nabla_{b_0}(\mathcal{R}^\nabla_{i}(b_{\sigma(1)},\cdots, b_{\sigma(k-j)}, \mathcal{R}^\nabla_{j}(b_{\sigma(k-j+1)}, \cdots,b_{\sigma(k-1)}, b_k), b_{k+1},\cdots, b_n)) \notag\\
 &=\sum_{i,j \geq 2, i+j=n+1}\sum_{k=j}^{n}\sum_{\sigma \in \sh(k-j,j-1)}\epsilon(\sigma)(-1)^{\abs{b_0} + \abs{b_{\sigma(1)}}+\cdots+\abs{b_{\sigma(k-j)}}} \notag \\
    &\qquad\qquad\mathcal{R}^\nabla_{i+1}(b_0,b_{\sigma(1)},\cdots,b_{\sigma(k-j)}, \mathcal{R}^\nabla_{j}(b_{\sigma(k-j+1)}, \cdots,b_{\sigma(k-1)}, b_k), b_{k+1},\cdots, b_n)\notag \\
 &\quad+ \sum_{i,j \geq 2, i+j=n+1}\sum_{k=j}^{n}\sum_{\sigma \in \sh(k-j,j-1)}\epsilon(\sigma)(-1)^{(\abs{b_{\sigma(1)}}+\cdots+ \abs{b_{\sigma(k-j)}})(\abs{b_0}+1)} \notag \\
    &\qquad\qquad\mathcal{R}^\nabla_{i}(b_{\sigma(1)},\cdots,b_{\sigma(k-j)}, \mathcal{R}^\nabla_{j+1} (b_0,b_{\sigma(k-j+1)}, \cdots,b_{\sigma(k-1)}, b_k), b_{k+1},\cdots, b_n).
\end{align}
Here in the last step we used the recursive definition of $\{\mathcal{R}^\nabla_i\}$.

Meanwhile, the second summand in Equation~\eqref{eq1} is
\begin{align}\label{eq3}
  &\quad[\mathcal{R}^\nabla_n,\mathcal{R}^\nabla_2(b_0,-)](b_1,b_2,\cdots,b_{n}) \notag \\
  &= \sum_{i=1}^n(-1)^{(\abs{b_0}+1)\ast_{i-1}}\mathcal{R}^\nabla_n(b_1,\cdots, \mathcal{R}^\nabla_2(b_0,b_i), \cdots,b_n) + (-1)^{\abs{b_0}}\mathcal{R}^\nabla_2(b_0,\mathcal{R}^\nabla_n(b_1,\cdots, b_{n})).
\end{align}
Substituting Equations~\eqref{eq2} and~\eqref{eq3} into Equation~\eqref{eq1}, we see that Equation~\eqref{casen} holds for the case $(n+1)$.
This proves that $\{\mathcal{R}^\nabla_k\}_{k\geq 1}$ satisfies Equation~\eqref{Leib} for all $n \geq 1$.
\ep
The Leibniz$_\infty[1] \A$-algebra $(\B,\{\mathcal{R}^\nabla_k\}_{k\geq 1})$ will be denoted by $\Kap^c(\delta)$. Here the superscript $c$ is to remind the reader that this Leibniz$_\infty[1] \A$-algebra is defined via a particular $\delta$-connection on $\B$.

\begin{Rem}
This method is originated from Kapranov's construction of $L_\infty$ algebra structure on the shifted tangent complex $\Omega_X^{0,\bullet-1}(T^{1,0}X)$ of a compact K\"{a}hler manifold $X$~\cite{Kap}. For this reason, we call $(\B,\{\mathcal{R}^\nabla_k\}_{k\geq 1})$  the Kapranov Leibniz$_\infty[1]$ $\A$-algebra.
\end{Rem}

\subsubsection{Functoriality}
Next, we show that the assignment of a Leibniz$_\infty[1]$ $\A$-algebra to each pair of dg derivation $\A \xrightarrow{\delta} \Omega$ and a $\delta$-connection $\nabla$ on $\B$ is functorial:
\begin{prop}\label{prop: functorial morphism}
  Let $\phi$ be a morphism from $\A \xrightarrow{\delta^\prime} \Omega^\prime$ to $\A \xrightarrow{\delta} \Omega$ in the category $\CATderivationA$ of dg derivations (see Definition~\ref{Def:morphismDerA}).
  Let $\B=\Omega^\vee$ and $\B'=(\Omega')^\vee$ be their dual dg $\A$-modules. For a $\delta$-connection $\nabla$ on $\B$ and a $\delta^\prime$-connection $\nabla^\prime$ on $\B^\prime$, there exists a morphism $f_\bullet = \{f_k\}_{k \geq 1}$ of Kapranov Leibniz$_\infty[1]$ $\A$-algebras from $(\B,\{\mathcal{R}^{\nabla}_k\}_{k\geq1})$ to $(\B^\prime, \{\mathcal{R}^{\nabla^\prime}_k\}_{k\geq1})$, whose first map is $f_1=\phi^\vee$.
  In other words, we have the following commutative diagram
  \bd
  \xymatrix{
    (\A \xrightarrow{\delta^\prime} \Omega^\prime) \ar[d]_-{\phi} \ar[r]^-{\Kap^{c}} & (\B^\prime,\{\mathcal{R}^{\nabla^\prime}_k\}_{k\geq1}) \\
    (\A \xrightarrow{\delta} \Omega) \ar[r]^-{\Kap^c} & (\B,\{\mathcal{R}^{\nabla}_k\}_{k \geq 1}). \ar[u]_-{\Kap^c(\phi) = f_\bullet}
  }
  \ed
\end{prop}
\bp
  Define a sequence of $\A$-multilinear maps $f_k: \B^{\otimes k} \rightarrow \B^\prime$ recursively by setting
 \begin{align}\label{phik}
 f_1(b_1) &= \phi^\vee(b_1), & f_{k+1}(b_0,\cdots,b_k) &= \nabla_{f_1(b_0)}^\prime f_k(b_1,\cdots,b_k) - f_k(\nabla_{b_0}(b_1,\cdots,b_k)), 
 \end{align}
 for all $k \geq 1$ and $b_i \in \B$. It is easy to verify that all the maps $\{f_k\}_{k\geq1}$ are $\A$-multilinear.
 Now we show that $\{f_k\}_{k\geq1}$ is a morphism of Leibniz$_\infty[1]$ $\A$-algebras from $(\B,\mathcal{R}_k^{\nabla})$ to $(\B^\prime, \mathcal{R}_k^{\nabla^\prime})$.

 We argue by induction: First of all, the $n=1$ case is obvious, since $f_{1} = \phi^\vee: \B \to \B^\prime$ is a morphism of dg modules.
 Now assume that Equation~\eqref{morphism equation} holds for some $n \geq 1$, i.e.,
 \begin{align}\label{casenprime}
  &\quad \sum_{k+p \leq n-1}\sum_{\sigma \in \sh(k,p)}\epsilon(\sigma)(-1)^{\dagger^{\sigma}_k}f_{n-p}(b_{\sigma(1)}, \cdots,b_{\sigma(k)}, \mathcal{R}^\nabla_{p+1}(b_{\sigma(k+1)}, \cdots,b_{\sigma(k+p)},b_{k+p+1}),\cdots,b_{n}) \notag \\
  &= \sum_{q \geq 1}\sum_{\substack{I^1 \cup \cdots \cup I^q = \mathbb{N}^{(n)} \\ I_1,\cdots,I_q \neq \emptyset \\ i_{\abs{I^1}}^1 < \cdots < i^q_{\abs{I^q}}}}\epsilon(I^1,\cdots,I^q) \mathcal{R}_q^{\nabla^\prime}(f_{\abs{I^1}}(b_{I^1}), \cdots,f_{\abs{I^q}}(b_{I^q})),
 \end{align}
 where $\dagger^{\sigma}_k = \sum_{i=1}^k\abs{b_{\sigma(i)}}$, $I^j = \{i^j_1 < \cdots < i^j_{\abs{I^j}}\} \subset \mathbb{N}^{(n)} = \{1,\cdots,n\}$, and $(b_{I^j}) = (b_{i^j_1},\cdots,b_{i^j_{\abs{I^j}}})$ for all $1\leq j \leq q$.

We proceed to show that Equation~\eqref{casenprime} holds for $(n+1)$ homogeneous inputs $\{b_0,\cdots,b_{n+1}\}$ .
For simplicity, we denote the left-hand side and the right-hand side of Equation~$\eqref{casenprime}$ by $\operatorname{LHS}(b_1,\cdots,b_n)$ and $\operatorname{RHS}(b_1,\cdots,b_n)$, respectively.

We write $\operatorname{LHS}(b_0,b_1,\cdots,b_n) = I_1 + I_2 + I_3$ as the sum of three parts, where
\begin{align*}
 I_1 &= f_{n+1}(\partial_{\A}b_0,b_1,\cdots,b_n) = \nabla^\prime_{\partial_\A f_1(b_0)}\circ f_n(b_1,\cdots,b_n) - \sum_{i=1}^n(-1)^{\ast_{i-1}(\abs{b_0}+1)} f_n(b_1,\cdots,\nabla_{\partial_\A b_0}(b_i),\cdots,b_n),\\
 I_2 &= \sum_{k+p =0}^{n-1}\sum_{\sigma\in \sh(k,p)}\epsilon(\sigma)(-1)^{\abs{b_0}+\dagger^{\sigma}_k} f_{n-p+1}(b_0,b_{\sigma(1)},\cdots,b_{\sigma(k)},\mathcal{R}^\nabla_{p+1} (b_{\sigma(k+1)},\cdots, b_{\sigma(k+p)},b_{k+p+1}),\cdots,b_n) \\
 &= (\nabla^\prime_{f_1(b_0)}\circ f_{n-p} - f_{n-p}\circ \nabla_{b_0})(b_{\sigma(1)}, \cdots, \mathcal{R}^\nabla_{p+1}(b_{\sigma(k+1)}, \cdots, b_{\sigma(k+p)},b_{k+p+1}),\cdots,b_n),
\end{align*}
by Equation~\eqref{phik}, and
\begin{align*}
  I_3 &= \sum_{k+p=0}^{n-1}\sum_{\sigma\in \sh(k,p)}\epsilon(\sigma)(-1)^{(\abs{b_0}+1)\dagger^{\sigma}_k }f_{n-p}(b_{\sigma(1)},\cdots,b_{\sigma(k)},\mathcal{R}^\nabla_{p+2}(b_0, b_{\sigma(k+1)},\cdots, b_{\sigma(k+p)},b_{k+p+1}),\cdots,b_n),\\
 &=\sum_{k+p=0}^{n-1}\sum_{\sigma \in \sh(k,p)} \epsilon(\sigma) (-1)^{(\abs{b_0}+1)\dagger^{\sigma}_k} 
 f_{n-p}(b_{\sigma(1)},\cdots, b_{\sigma(k)},-[\mathcal{R}^\nabla_{p+1}, \nabla_{b_0}](b_{\sigma(k+1)},\cdots, b_{\sigma(k+p)},b_{k+p+1}),\cdots,b_n) \\
  &\quad\quad +\sum_{k=0}^{n-1}(-1)^{\ast_{k}(\abs{b_0}+1)} f_n(b_1,\cdots,\nabla_{\partial_\A b_0}(b_{k+1}),\cdots,b_n),
\end{align*}
by Equations~\eqref{Atiyah cocycle} and~\eqref{Rnabla}.
Summing them up, we have
 \begin{align}\label{LHS}
  \operatorname{LHS}(b_0,b_1,\cdots,b_n)
  &= \nabla^\prime_{\partial_\A(f_1(b_0))} f_n(b_1,\cdots,b_n) -\sum_{i=1}^n(-1)^{\abs{b_0}\ast_{i-1}}\operatorname{LHS}(b_1, \cdots,\nabla_{b_0}b_i,\cdots,b_n) \notag \\
  &\quad + \nabla^\prime_{f_1(b_0)} \operatorname{LHS}(b_1,\cdots,b_n).
 \end{align}

Meanwhile,
\begin{align}
  &\quad\operatorname{RHS}(b_0,b_1,\cdots,b_n) \notag\\ \nonumber
  &= \sum_{q=1}^{n}\sum_{\substack{I^1\cup\cdots\cup I^q = \mathbb{N}^{(n)} \\ I_1,\cdots,I_q \neq \emptyset \\ i_{\abs{I^1}}^1 < \cdots < i^q_{\abs{I^q}}}} \epsilon(I^1,\cdots,I^q) (\mathcal{R}^{\nabla^\prime}_{q+1}(b_0, f_{\abs{I^1}}(b_{I^1}),\cdots,f_{\abs{I^q}}(b_{I^q})) \\\nonumber
  &\qquad+ \mathcal{R}^{\nabla^\prime}_{q} (f_{\abs{I^1}+1}(b_0,b_{I^1}), \cdots, f_{\abs{I^q}}(b_{I^q})) + \cdots + (-1)^{\abs{b_0}(\abs{b_{I^1}}+\cdots+ \abs{b_{I^{q-1}}})} \mathcal{R}^{\nabla^\prime}_{q} (f_{\abs{I^1}}(b_{I^1}),\cdots, f_{\abs{I^q}+1}(b_0,b_{I^q}))) \notag\\
  &\qquad\qquad\qquad\text{by Equations~\eqref{Atiyah cocycle},\eqref{Rnabla},\eqref{phik}}  \notag \\  \label{RHS}
  &=\nabla^\prime_{f_1(b_0)} \operatorname{RHS}(b_1,\cdots,b_n) + \nabla^\prime_{\partial_\A(f_1(b_0))}f_n(b_1,\cdots,b_n)
  -\sum_{i=1}^n(-1)^{\abs{b_0}\ast_{i-1}}\operatorname{RHS}(b_1,\cdots, \nabla_{b_0}b_i, \cdots, b_n).
\end{align}
Applying the induction assumption to Equations~\eqref{LHS} and~\eqref{RHS}, we see that Equation~\eqref{casenprime} holds for all $(n+1)$ entries.
Thus Equation~\eqref{casenprime} holds for all $n \geq 1$. This proves that $f = \{f_k\}_{k\geq1}$ is a morphism of Leibniz$_\infty[1]$ $\A$-algebras.
\ep

As a consequence, the Kapranov's construction defines a contravariant functor
  \bd
   \Kap^c:  \CATderivationA \to \CATLeibnizoneA
  \ed
from the category $\CATderivationA$ of dg derivations of $\A$ to the category $\CATLeibnizoneA$ of Leibniz$_\infty[1]$-algebras.

\begin{Rem}\label{Rmk:whyoverA}
The reason that we restrict to work in the category $\CATLeibnizoneA$ of Leibniz$_\infty[1]$ $\A$-algebras is as follows:
If we treat $(\B,\{\mathcal{R}^\nabla_k\}_{k\geq 1})$ merely as a Leibniz$_\infty[1]$ algebra over $\k$, it is always isomorphic to the trivial one $(B,\{\partial_\A,0,0,\cdots\})$ (all higher brackets are zero). In fact,
 one can build a sequences of degree $0$  maps
\begin{align*}
  \phi_k:~\; &\B^{\otimes k} \rightarrow \B,\quad k \geq 1,
\end{align*}
where $\phi_1 = \id_\B$, and $\{\phi_{k+1}\}_{k\geq1}$ are defined recursively by
\begin{align*}
 \phi_{k+1}(b_0,\cdots,b_k) &=   {\nabla}_{b_0} \circ \phi_k (b_1,\cdots,b_k), \;\;\;\forall b_i \in \B.
\end{align*}
The set $\{\phi_k:~ \B^{\otimes k} \rightarrow \B\}_{k\geq1}$  defines an isomorphism of Leibniz$_\infty[1]$ algebras from $(\B,\{\partial_\A,0,0,\cdots\})$ to $(\B,\{\mathcal{R}^{\nabla}_k\}_{k\geq1})$ in the category of Leibniz$_\infty[1]$ algebras over $\k$. The proof is similar to that of Proposition~\ref{prop: functorial morphism}.
However, the maps $\{\phi_k\}_{k\geq 2}$ are not $\A$-multilinear.
\end{Rem}

Next, we stress the independence from the choice of connections in the definition of Kapranov functors. For a dg derivation $\A \xrightarrow{\delta} \Omega$, suppose that we have another $\delta$-connection $\widetilde{\nabla}$ on $\B = \Omega^\vee$. Denote the corresponding Kapranov Leibniz$_\infty[1]$ $\A$-algebra by $\Kap^{\tilde{c}}(\delta) = (\B,\mathcal{R}_k^{\widetilde{\nabla}})$.
By Proposition~\ref{prop: functorial morphism}, there exists an isomorphism $g^{\nabla,\widetilde{\nabla}}_\bullet: \Kap^c(\delta) \to \Kap^{\tilde{c}}(\delta)$ of Leibniz$_\infty[1]$ $\A$-algebras, where $g^{\nabla,\widetilde{\nabla}}_1 = \id_\B$, and $\{g^{\nabla,\widetilde{\nabla}}_{k+1}\}_{k\geq1}$ are defined recursively as follows:
\begin{align*}
 g^{\nabla,\widetilde{\nabla}}_{k+1}(b_0,\cdots,b_k) &= (\widetilde{\nabla}_{b_0} \circ g^{\nabla,\widetilde{\nabla}}_k - g^{\nabla,\widetilde{\nabla}}_k \circ \nabla_{b_0})(b_1,\cdots,b_k), \;\;\;\forall b_i \in \B.
\end{align*}
Moreover, via a straightforward verification, we have
\begin{lem}\label{independence of connection}
 There exists a natural equivalence between Kapranov functors $\Kap^c$ and $\Kap^{\tilde{c}}$ with respect to different connections. In other words, for any morphism $\phi: (\A \xrightarrow{\delta^\prime} \Omega^\prime) \to (\A \xrightarrow{\delta} \Omega)$ of dg derivations of $\A$, we have the following commutative diagram
 \bd
  \xymatrix{
    \Kap^c(\delta) \ar[d]_-{\Kap^c(\phi)} \ar[r]^-{g^{\nabla,\widetilde{\nabla}}_\bullet}  & \Kap^{\tilde{c}}(\delta) \ar[d]^-{\Kap^{\tilde{c}}(\phi)} \\
    \Kap^c(\delta^\prime) \ar[r]^-{g^{\nabla^\prime,\widetilde{\nabla^\prime}}} & \Kap^{\tilde{c}}(\delta^\prime).
  }
 \ed
\end{lem}
By this natural equivalence, we are allowed to drop the superscript $c$ to obtain the following
\begin{Thm}\label{Thm: Kapranov functor}
  The Kapranov's construction defines a contravariant functor
  \bd
   \Kap:  \CATderivationA \to \CATLeibnizoneA
  \ed
  from the category $\CATderivationA$ of dg derivations of $\A$ to the category $\CATLeibnizoneA$ of Leibniz$_\infty[1]$-algebras.
\end{Thm}
\begin{Rem}
By the universal property, the K\"{a}hler differential $\A \xrightarrow{d_{dR}} \Omega_{\A\mid\k}^1$ is the initial object in the category $\CATderivationA$ of dg derivations. Thus the corresponding Kapranov Leibniz$_\infty[1]$ $\A$-algebra on the tangent complex $T_{\A\mid\k} = (\Omega_{\A\mid\k}^1)^\vee$ of $\A$ is the final object of the subcategory in $\CATLeibnizoneA$ consisting of Kapranov Leibniz$_\infty[1]$ $\A$-algebras arising from dg derivations of $\A$.
\end{Rem}

Let $\A \xrightarrow{\delta} \Omega$ be a dg derivation of $\A$ and $\E$ a dg $\A$-module. By a similar argument, $\E$ carries a Leibniz$_\infty[1]$ $\A$-module structure over $\Kap(\delta)$. Moreover, we have
\begin{Thm}\label{Thm: module}
Given a dg derivation $\A \xrightarrow{\delta} \Omega$ of $\A$, there exists a functor from 
the category $\DG\A$ of dg $\A$-modules to the category of Leibniz$_\infty[1]$ $\A$-modules over $\Kap(\delta)$.
\end{Thm}

\subsubsection{Leibniz algebra structures}
Let $(V,\{\lambda_k\}_{k\geq1})$ be a Leibniz$_\infty[1]$ $\A$-algebra as in Definition \ref{Def:LeibnizinfinityA}. Then $(V,\lambda_1=\partial_\A)$ is a dg $\A$-module. Its cohomology $H^\bullet(V)$ is called the tangent cohomology of the Leibniz$_\infty[1]$ $\A$-algebra $(V,\{\lambda_k\}_{k\geq1})$.
According to~\cite{CSX}*{Proposition 3.10}, the (degree $(-1)$ shifted) tangent cohomology $H^\bullet(V[-1])$ is a Leibniz algebra (over $\k$), when equipped with the bracket
$$
 \check{\lambda}_2:~H^\bullet(V[-1]) \times H^\bullet(V[-1]) \to H^\bullet(V[-1])
$$
$$
\check{\lambda}_2([x],[y]):=(-1)^{\abs{x}}[\lambda_2(x,y)],
$$
where $x,y\in V$ are $\lambda_1$-closed.

In a similar fashion,   if $(W,\{\mu_k\}_{k\geq1})$ is a $(V,\{\lambda_k\}_{k\geq1})$ $\A$-module as in Definition \ref{Def:LeibnizinfinityAmodule}, then   $(W,\mu_1=\partial_\A)$ is also a dg $\A$-module. The cohomology $H^\bullet(W)$  is a Leibniz module over the aforesaid Leibniz algebra $ H^\bullet(V[-1])$ (both over $\k$), when equipped with the action
$$ \check{\mu}_2:~H^\bullet(V[-1]) \times H^\bullet(W) \to H^\bullet(W)$$
$$\check{\mu}_2([x],[w]):=(-1)^{\abs{x}}[\mu_2(x,w)],$$
where $x\in V,~ w\in W$ are, respectively, $\lambda_1$- and $\mu_1$-closed elements.

As a consequence of Theorem~\ref{Thm: Functoriality}, Theorems~\ref{Thm: leibniz infty} and~\ref{Thm: module}, we have the following
\begin{Cor}\label{maincoro}
  Let $\phi$ be a morphism of dg derivations from $\A \xrightarrow{\delta^\prime} \Omega^\prime$ to $\A \xrightarrow{\delta} \Omega$ and let $\B$ and $\B^\prime$ be the dual dg $\A$-modules of $\Omega$ and $\Omega^\prime$, respectively.
  \begin{compactenum}
    \item The (degree $(-1)$ shifted) cohomology space $H^\bullet(\A,\B[-1])$ is a Leibniz algebra, whose bracket ${\Bigl[-,-\Bigr]}_\B$ is induced by the $\delta$-twisted Atiyah class of $\B$:
        $$
        \Bigl[[b_1 ],[b_2 ]\Bigr]_\B = (-1)^{\abs{b_1}}\At^\delta_\B ([b_1],[b_2]),
        $$
        where $b_1,b_2 \in \B$ are $\partial_\A$-closed elements. Moreover, $\phi^\vee:~\B\to \B^\prime$ induces a morphism of Leibniz algebras, i.e.,
        \bd
         \Bigl[\phi^\vee(b_1 ),\phi^\vee[b_2]\Bigr]_{\B^\prime} = \phi^\vee(\Bigl[[b_1 ],[b_2 ]\Bigr]_\B).
        \ed
    \item For any dg $\A$-module $\E$, there exists a representation of $H^\bullet(\A,\B[-1])$ on the cohomology space $H^\bullet(\A,\E)$, with the action map $-\triangleright-$ induced by the $\delta$-twisted Atiyah class of $\E$:
        $$
        [b] \triangleright [e]= (-1)^{\abs{b}}\At^\delta_\E ([b],[e]),
        $$
        where $b\in\B$, $e\in\E$ are both $\partial_\A$-closed elements. Moreover, this assignment is functorial, i.e., for each dg $\A$-module morphism $\lambda: \E  \to \F$ (of degree $0$),
        \bd
          [b] \triangleright \lambda(e) = \lambda ([b] \triangleright [e]).
        \ed
   \end{compactenum}
\end{Cor}

\begin{Rem}
 According to \cite{CSX}*{Theorem 3.4}, the Atiyah class of a Lie pair  $(L,A)$ induces a Lie algebra structure on the cohomology $H^\bullet_{\CE}(A, L/A[-1])$. A similar result holds for $L_\infty$ algebra pairs~\cite{CLX}. However, it is not the case in general (see an example below). It is natural to ask when the Leibniz algebra structure in Corollary~\ref{maincoro} could be refined to a Lie algebra structure. We will investigate this question somewhere else.
\end{Rem}

\begin{Ex}
Let $\mathcal{LM}$ be the category of linear maps~\cite{LP}. A Lie algebra object in $\mathcal{LM}$ is a triple $E \xrightarrow{\psi} \g$, where $\g$ is a Lie algebra, $E$ is a left $\g$-module, and $\psi$ is a $\g$-equivariant linear map.
Consider the cdga $\A=C^\bullet(\g) = (\wedge^\bullet\g^\vee,d_{\CE})$ and dg $C^\bullet(\g)$-module $\Omega=C^{\bullet}(\g,E^\vee[-1]) = (\wedge^{\bullet}\g^\vee \otimes E^\vee[-1], d_{\CE} )$, i.e., the Chevalley-Eilenberg cochain complex of the dual $\g$-module $E^\vee[-1]$.
The $\g$-equivariant map $E \xrightarrow{\psi} \g$ gives rise to a dg derivation of $C^\bullet(\g) $:
\bd
  C^\bullet(\g) \xrightarrow{\delta = \psi^\vee} C^{\bullet}(\g,E^\vee[-1]).
\ed
The dual module of $\Omega=C^{\bullet}(\g,E^\vee[-1]) $ is $\B=C^{\bullet}(\g,E[1])$.
One can take the trivial $\delta$-connection on $\B$:
\bd
 \nabla: \B=C^\bullet(\g,E[1]) \rightarrow \Omega\otimes_\A \B=C^{\bullet}(\g,E^\vee[-1] \otimes E[1]),
\ed
defined by
\bd
 \nabla(\omega \otimes e) = \delta(\omega) \otimes e,\;\;\forall \omega \in \wedge^\bullet\g^\vee, e \in E.
\ed
By Equation~\eqref{Atiyah cocycle}, the associated Atiyah cocycle is a degree $1$ element $\At^\nabla_\B \in E^\vee[-1]\otimes E^\vee[-1] \otimes E[1]$  specified by
\bd
 \mathcal{R}^\nabla_2=\At^\nabla_\B(e_1,e_2) = -\psi(e_1)e_2 ,\;\;\;\forall e_1,e_2 \in E .
\ed

It can be easily seen that higher structures $\mathcal{R}^\nabla_j=0$ for all $j\geq 3$. Hence,  the Kapranov Leibniz$_\infty[1]$ $C^\bullet(\g)$-algebra $\B=C^\bullet(\g,E[1])$ is simply a dg Leibniz$[1]$ algebra in this case, or equivalently, $\B[-1]=C^\bullet(\g,E)$ is a dg Leibniz  algebra. In particular, the subspace $E$ is a Leibniz algebra, recovering the result in~\cite{LP}.

By Corollary~\ref{maincoro}, there is a Leibniz algebra structure on the graded vector space $H^\bullet(\A,\B[-1])=H_{\CE}^\bullet(\g,E)$, whose bracket is given by
\bd
   \Bigl[[e_1],[e_2]\Bigr] = (-1)^{\abs{e_1}+1}[\At^\nabla_\B(e_1 ,e_2 ) ] =\pm [\psi(e_1)e_2],
\ed
for all $d_{\CE}$-closed elements $e_1,e_2 \in C^\bullet(\g,E )$. Here the last term $\psi(-)(-):~C^\bullet(\g,E)\times C^\bullet(\g,E )\to C^\bullet(\g,E )$ is a $(\wedge^\bullet\g^\vee)$-bilinear map naturally extended from $\psi(-)(-):~E\times E\to E$.
In general, the Leibniz structure on $(H_{\CE}^\bullet(\g,E),[-,-])$ is not skewsymmetric.
\end{Ex}

\subsubsection{Homotopic invariance}
In this section, we prove that the isomorphism class of Kapranov Leibniz$_\infty[1]$ $\A$-algebras arising from dg derivations only depends on their homotopy classes.
\begin{prop}\label{Thm: homotopy invariance}
Let $\delta \sim \delta^\prime$ be homotopic $\Omega$-valued dg derivations of $\A$.Then there exists an isomorphism $\{g_k\}_{k\geq1}$ sending the Kapranov Leibniz$_\infty[1]$ $\A$-algebra $\Kap(\delta^\prime) = (\B,\{\mathcal{R}^{\nabla^\prime}_k\}_{k\geq1})$ (with respect to a $\delta^\prime$-connection $\nabla^\prime$) to $\Kap(\delta) = (\B,\{\mathcal{R}^{\nabla}_k\}_{k\geq1})$ (with respect to a $\delta$-connection $\nabla$).
\end{prop}
\bp
By assumption, there exists a degree $(-1)$ $\Omega$-valued derivation $h:\A  \rightarrow \Omega$ of $\A$ such that
$$
\delta^\prime = \delta + [\partial_\A,h] = \delta + \partial_\A  \circ h + h \circ d_\A.
$$
We choose an $h$-connection on $\B$, i.e. a degree $(-1)$ linear map
\bd
\widehat{\nabla}: \B \rightarrow \Omega \otimes_\A \B
\ed
satisfying
$$
\widehat{\nabla}(ab) = h(a)\otimes b + (-1)^{\abs{a}}a \widehat{\nabla}(b),\;\; \forall a \in \A, b \in \B.
$$
For each $\delta$-connection $\nabla$ on $\B$, it can be easily verified that $\nabla^{\prime\prime} := [\partial_\A,\widehat{\nabla}]$ is a $[\partial_\A,h]$-connection on $\B$, and thus
\bd
 \nabla^\prime = \nabla + \nabla^{\prime\prime} = \nabla + [\partial_\A,\widehat{\nabla}]: \B \rightarrow \Omega \otimes_\A \B
\ed
is a $\delta^\prime$-connection on $\B$. It follows that
\bd
 \mathcal{R}^{\nabla^\prime}_2 = [\nabla^\prime,\partial_\A] = [\nabla + [\partial_\A, \widehat{\nabla}],\partial_\A] = [\nabla, \partial_\A] = \mathcal{R}_2^{\nabla}.
\ed
Define a family of $\A$-multilinear maps $g_k: \B^{\otimes k} \rightarrow \B$ inductively by setting $g_1 = \id_\B, g_2 = 0$, and
\begin{align}\label{gk}
 g_{k+1}(b_0,\cdots,b_k) &= (-1)^{\abs{b_0}}\sum_{p=2}^k \sum_{\substack{I^1\cup\cdots\cup I^p = \mathbb{N}^{(k)} \\ I_1,\cdots,I_p \neq \emptyset \\ i_{\abs{I^1}}^1 < \cdots < i^p_{\abs{I^p}}}} \epsilon(I^1,\cdots,I^p) [\widehat{\nabla}_{b_0},R^{\nabla}_{p}](g_{\abs{I^1}}(b_{I^1}),\cdots, g_{\abs{I^p}}(b_{I^p})) \notag \\
 &\quad + [\nabla^\prime_{b_0}, g_k](b_1,\cdots,b_k),
\end{align}
for all $k\geq 2$.
It follows from a straightforward inductive argument that $\{g_k\}_{k\geq 1}$ is a morphism of Leibniz$_\infty[1]$ $\A$-algebras from $(\B,\mathcal{R}_n^{\nabla^\prime})$ to $(\B,\mathcal{R}_n^\nabla)$.
\ep

\begin{Rem}
  Although the Kapranov functor $\Kap$ maps homotopic derivations to isomorphic Leibniz$_\infty[1]$ $\A$-algebra, it does not reduce to a functor from the category consisting of homology classes of dg derivations of $\A$ to the category $\CATLeibnizoneA$ of Leibniz$_\infty[1]$ $\A$-algebras.
\end{Rem}

\subsection{Applications}
We first consider a Lie pair $(L,A)$, and let $B=L/A$ be the Bott $A$-module. In the introduction, we explained that for each splitting $j: B \to L$ of the short exact sequence~\eqref{SES} and for any $L$-connection $\nabla$ on $B$ extending the Bott $A$-module structure, there associates a Leibniz$_\infty[1]$ algebra structure $\{\lambda_k\}_{k\geq1}$ on the graded $\k$-vector space $\Gamma(\wedge^\bullet A^\vee \otimes B)$. As all $\{\lambda_k\}_{k\geq2}$ are $\Omega_A^\bullet$-multilinear, it is a
Leibniz$_\infty[1]$ $\Omega_A^\bullet$-algebra.

Recall that we have a $\Omega_A^\bullet(B^\vee)$-valued dg derivation $\delta_j$ of the cdga $\Omega_A^\bullet$ as in Equation~\eqref{Eqt:deltaofLA}. By Proposition~\ref{prop:Lie pairs}, the Atiyah cocycle $\alpha_B^\nabla$ of the Lie pair coincides with the Atiyah cocycle $\At^{\nabla^{\delta_j}}_\B$ of the dg $\Omega_A^\bullet$-module $\B:= \Omega_A^\bullet(B)$ with respect to a $\delta_j$-connection $\nabla^{\delta_j}$ as in Equation~\eqref{Eq:Landdeltaconnection}. Comparing definitions of $\{\lambda_k\}_{k\geq3}$ in the introduction and $\{\mathcal{R}^{\nabla^{\delta_j}}_k\}$ as in Equation~\eqref{Rnabla}, we see that the two Leibniz$_\infty[1]$ $\Omega_A^\bullet$-algebras $(\B,\{\lambda_k\}_{k\geq1})$ and $(\B,\mathcal{R}^{\nabla^{\delta_j}}_k)$ are exactly the same.

Applying Proposition~\ref{prop: splitting}, Theorem~\ref{Thm: Kapranov functor}, Theorem~\ref{Thm: module}, and Proposition~\ref{Thm: homotopy invariance}, we have the following

\begin{Thm}\label{Thm:Liepair}
  Let $(L,A)$ be a Lie  pair  over a smooth manifold $M$. Then the Leibniz$_\infty[1]$ algebra structure constructed in~\cite{CSX}*{Theorem 3.13} on the graded $\k$-vector space $\Gamma(\wedge^\bullet A^\vee \otimes L/A)$ is unique up to isomorphisms in the category of Leibniz$_\infty[1]$  $\Omega_A^\bullet$-algebras.

  Moreover, if $(E,\partial_A^E)$ is an $A$-module, then the representation of the above Leibniz$_\infty[1]$ algebra on the graded $\k$-vector space $\Gamma(\wedge^\bullet A^\vee \otimes E)$   is also unique up to isomorphisms in the category of Leibniz$_\infty[1]$  $\Omega_A^\bullet$-modules.
\end{Thm}

Finally, we consider another interesting application:
Let $X$ be a complex manifold and $\A = (\Omega_X^{0,\bullet}, \bar{\partial})$ its Dolbeault dg algebra. Let $\Omega = (\Omega_X^{0,\bullet}(T^{1,0}X),\bar{\partial})$ be the dg $\A$-module generated by the smooth section space $\Gamma(T^{1,0}X)$ of the holomorphic tangent bundle $T^{1,0}X$.
Note that each holomorphic bivector field $\pi \in \Gamma(\wedge^2 T^{1,0}X)$ determines an $\Omega$-valued dg derivation of $\A$, denoted by $\delta_\pi$, which is the composition
\bd
  \A \xrightarrow{\partial} \Omega_X^{1,\bullet} = \Omega_X^{0,\bullet}((T^{1,0}X)^\vee) \xrightarrow{\pi^\sharp} \Omega.
\ed
Here $\pi^\sharp$ is the contraction along $\pi$ from $(T^{1,0}X)^\vee$ to $T^{1,0}X$.

In fact, $\pi^\sharp$ is a morphism of dg derivations of $\A$ (from $\A \xrightarrow{\partial} \Omega_X^{1,\bullet}$ to $\A \xrightarrow{\delta_\pi} \Omega$). It sends the Atiyah class $\alpha_E \in H^1(X,(T^{1,0}X)^\vee\otimes\End(E))$ of any holomorphic vector bundle $E$ to the $\delta_\pi$-twisted Atiyah class $\At^{\delta_\pi}_\E \in H^1(X,T^{1,0}X \otimes \End(E))$ of the associated dg $\A$-module $\E = \Omega_X^{0,\bullet}(E)$.

By Proposition~\ref{prop: Atiyah vanish}, the $\delta_\pi$-twisted Atiyah class $\At^{\delta_\pi}_\E$ measures the existence of holomorphic $\delta_\pi$-connections on $E$. In particular, if $\pi$ a holomorphic Poisson bivector field, then $(T^{1,0}X)^\vee$ is a holomorphic Lie algebroid~\cite{LGSX}, and $\At^{\delta_\pi}_\E$ measures the existence of holomorphic $(T^{1,0}X)^\vee$-connections on $E$.

Applying Theorem~\ref{Thm: Kapranov functor}, we have the following
\begin{Thm}
  Let $X$ be a complex manifold,  $\pi$ a holomorphic bivector field. Then,
  \begin{itemize}
  \item Both $\Omega_X^{0,\bullet}(T^{1,0}X)$ and $\Omega_X^{0,\bullet}((T^{1,0}X)^\vee)$ carry canonical Kapranov Leibniz$_\infty[1]$ $\Omega_X^{0,\bullet}$-algebra structures;
  \item There is a morphism of Leibniz$_\infty[1]$ $\Omega_X^{0,\bullet}$-algebras $\{f_k\}_{k\geq1}: \Omega_X^{0,\bullet}((T^{1,0}X)^\vee) \to \Omega_X^{0,\bullet}(T^{1,0}X)$ such that $f_1 = \pi^\sharp$.
  \end{itemize}
\end{Thm}

\section{Open questions and remarks}\label{Sec: final section}
In this note, we assume that each dg $\A$-module $\E$ is projective in order that connections exist on $\E$. In the non-projective case, one can follow Calaque-Van den Bergh's approach \cites{CV} to define the Atiyah class of $\E$ (which coincides with the Atiyah class of $\E$ in Definition~\ref{Def Costello} when $\E$ admits connections)--- The first step is to construct a short exact sequence, called the jet sequence, of dg $\A$-modules:
 $$
 \xymatrix@C=0.5cm{
   0 \ar[r] & \Omega_{\A\mid\k}^1\otimes_\A \E \ar[rr]  && \mathfrak{J}\E \ar[rr]  && \E \ar[r] & 0 }.
 $$
The Atiyah class of $\E$ is then defined to be the extension class of the above jet sequence. We would like to follow this approach to study twisted Atiyah classes of some cases when connections do not exist (singular foliations considered in~\cite{LGLS} for example).

Note that Kapranov's original construction on $\Omega^{0,\bullet-1}_X(T_X)$ of a K\"{a}hler manifold $X$ is  an $L_\infty$ algebra, whereas Chen, Sti\'{e}non and Xu's construction of $\Gamma(\wedge^\bullet A^\vee \otimes B)$ is a Leibniz$_\infty[1]$ algebra. In fact, this is due to the existence of Chern connection on $T_X$ which enjoys special properties (see \cite{CSX}*{Section 3.4.4}).
Meanwhile, when $\A = C^\infty(\M)$ is the cdga of functions of a smooth dg manifold $\M$. According to~\cite{MSX}, the tangent complex $T_{\A\mid\k} = \Gamma(T_\M)$ admits an $L_\infty[1]$ algebra structure (by a construction different from the Kapranov's construction we discussed).
Moreover, Laurent-Gengoux, Sti\'{e}non and Xu~\cite{LaurentSX} have proved that for each Lie pair $(L,A)$, there exists a canonical $L_\infty[1]$ algebra structure on the graded $\k$-vector space $\Gamma(\wedge^\bullet A^\vee \otimes L/A)$ (which is different from the Chen-Sti\'{e}non-Xu's construction in \cite{CSX}). It is natural to ask how to tweak the Kapranov Leibniz$_\infty[1]$ algebra of general dg derivations so as to produce an $L_\infty[1]$ algebra rather than a mere Leibniz$_\infty[1]$ algebra.

According to the perturbation lemmas proved by Huebschmann~\cites{HueLie,HueshLie}, many $L_\infty$ algebras arise from dg Lie algebras or $L_\infty$ algebras by homological perturbation theory. It is interesting to investigate whether similar perturbation lemma holds for Leibniz$_\infty[1]$ algebras. Moreover, if this is the case, then it is natural to ask for which kind of dg derivations of a cdga $\A$, the associated Kapranov Leibniz$_\infty[1]$ $\A$-algebra results from some perturbation.

These questions will be investigated somewhere else.

We would also like to mention other works that are related to the present paper:
Batakidis and Voglaire~\cite{BV} showed how Atiyah classes of Lie pairs~\cite{CSX} and of dg Lie algebroids~\cite{MSX} give rises to Atiyah classes of dDG algebras~\cite{CV}.
Bordemann~\cite{Bordemann} studied the Atiyah class as the obstruction to the existence of invariant connections on homogeneous spaces.
Hennion~\cite{Hen} generalized Kapranov's construction to algebraic derived stack: There exists a Lie algebra structure on the shifted tangent complex $\mathbb{T}_X[-1]$ of a derived Artin stack $X$ locally of finite presentation. Moreover, given a perfect module $E$ over $X$, there exists a representation of the aforesaid Lie algebra on $E$ induced by the Atiyah class of $E$.

\end{document}